\documentclass{amsart}

\usepackage{amsmath, amscd, amsthm, amssymb, graphics, mathrsfs,url}

\theoremstyle{plain}
\newtheorem{thm}{Theorem}[subsection]
\newtheorem{lem}[thm]{Lemma}
\newtheorem{cor}[thm]{Corollary}
\newtheorem{prop}[thm]{Proposition}

\theoremstyle{definition}
\newtheorem{defn}[thm]{Definition}

\theoremstyle{remark}
\newtheorem{rem}[thm]{Remark}
\newtheorem{ex}[thm]{Example}

\newtheorem{exampleit}[thm]{Example}

\newenvironment{example}{\begin{exampleit}\rm}{\end{exampleit}}

\newcommand{\AAA}{{\mathbb{A} }}

\newcommand{\CC}{{\mathbb{C} }}

\newcommand{\ZZ}{{\mathbb{Z} }}
\newcommand{\PP}{ {\mathbb{P} } }
\newcommand{\QQ}{{\mathbb{Q} }}

\def\a{\alpha}

\def\im{\mathit{im}}

\newcommand{\liets}{{\bf t}^*}
\newcommand{\Proj}{{\rm Proj}}
\newcommand{\Spec}{{\rm Spec}}
\newcommand{\diag}{{\rm diag}}

\title{Towards Non-Reductive Geometric Invariant Theory}

\author{Brent Doran}

\address{School of Mathematics \\ Institute for Advanced Study \\
Princeton, NJ 08540, USA}

\email{doranb@math.ias.edu}

\author{Frances Kirwan}
\address{Mathematical Institute \\ University of Oxford \\ Oxford, OX1 3LB, UK}

\email{kirwan@maths.ox.ac.uk}

\date{\today}

\dedicatory{To Bob MacPherson on his 60th birthday}

\begin{document}

\begin{abstract}
We study linear actions of algebraic groups on smooth projective
varieties $X$. A guiding goal for us is to understand the cohomology
of ``quotients" under such actions, by generalizing (from reductive
to non-reductive group actions) existing methods involving Mumford's
geometric invariant theory (GIT). We concentrate on actions of
unipotent groups H, and define sets of stable points $X^s$ and
semistable points $X^{ss}$, often explicitly computable via the
methods of reductive GIT, which reduce to the standard definitions
due to Mumford in the case of reductive actions. We compare these
with definitions in the literature. Results include (1) a geometric
criterion determining whether or not a ring of invariants is
finitely generated, (2) the existence of a geometric quotient of
$X^s$, and (3) the existence of a canonical ``enveloping quotient"
variety of $X^{ss}$, denoted $X/\!/H$, which (4) has a projective
completion given by a reductive GIT quotient and (5) is itself
projective and isomorphic to $\Proj(k[X]^H)$ when $k[X]^H$ is
finitely generated.
\end{abstract}

\maketitle

\begin{footnotesize}
\tableofcontents
\end{footnotesize}

\section{Introduction}

Geometric invariant theory (GIT) is a powerful theoretical and
computational tool for the study of reductive algebraic group
actions.  On the theory side, it provides a good notion of a
quotient of an affine or projective variety; many key moduli spaces
admit a description as a GIT quotient. More computationally,
information about the geometry of projective GIT quotients, in
particular about their cohomology (indeed, algebraic K-theory,
motivic cohomology or any oriented cohomology theory
\cite{ADK1,Kir1,ES,ChaiNeeman,EdGr}), can often be extracted using
not much more than the theory of weights of linear representations.
As an added bonus, the ``stable'' structures identified by GIT tend
to have very nice (differential) geometric properties \cite{Thomas}.
A drawback is that GIT requires the group which acts to be
reductive, whereas many interesting problems in moduli theory,
affine geometry, and classical invariant theory are rooted in
non-reductive actions. This paper addresses the question of how to
develop an effective version of GIT for general affine algebraic
group actions, including some of the difficulties, context, and
motivating problems, with the eventual goal of computing the
cohomology of quotients of nonsingular projective varieties by such
actions.

Two key properties of actions of a reductive group $G$ on an affine
variety $X$ are crucial to GIT.  Firstly, by a theorem of Nagata,
the ring of invariants $k[X]^G$ is finitely generated. Secondly,
given any two disjoint closed $G$-invariant subvarieties of $X$,
there exists an invariant function which separates them.  These lead
to the existence of a ``good categorical quotient'' $\pi: X
\rightarrow X/\!/G = \Spec(k[X]^G)$, along with a distinguished open
subset $X^s$ of {\em stable} orbits such that $\pi(X^s) \subseteq
X/\!/G$ is a ``geometric quotient'', that is, an orbit space with
nice properties. Quotients of more general varieties $X$ (equipped
with linearizations of the group actions) are constructed by
patching together suitable open affine pieces to get a categorical
quotient $X/\!/G$ of the {\em semistable} subset $X^{ss}$ of $X$,
with an open subset $X^s/G \subseteq X/\!/G$ giving a geometric
quotient of $X^s$.

Less essential, but still very useful, is a third property of
reductive actions, that every invariant extends to an invariant of
the ambient affine (or projective) space. This effectively reduces
GIT to the study of representations of reductive groups on the
affine (or projective) space itself. Fourthly and finally, the
quotients $X/\!/G$ of affine or projective varieties are again
affine or projective, respectively, with the quotient maps $X^{ss}
\to X/\!/G$ surjective.

Each of these properties fails for non-reductive groups.  Most
famously, examples of non-finitely generated rings of invariants
(counter-examples to Hilbert's Fourteenth Problem) were first
discovered by Nagata \cite{Nagata}. One might hope, however, that a
good generalization of GIT to non-reductive group actions would, as
in the reductive theory, satisfy the following properties:

\begin{itemize}
\item It would depend only on the data of a
$G$-linearization of $X$; that is, a $G$-equivariant embedding of
$X$ into an affine or projective space equipped with a linear action
of $G$, together with a lift of the action of $G$ to the
(homogeneous) coordinate ring $k[X]$.
\item It would use the invariants in $k[X]$
 to separate as many closed orbits
as possible.
\item It would provide notions of stable and semistable subsets
$X^s$ and $X^{ss}$ with a canonical $G$-invariant morphism $X^{ss}
\to X/\!/G$ (ideally a categorical quotient) restricting to a
geometric quotient of $X^s$.
\item It would have a good `change of groups' formalism relating
 $X/\!/H$ and $(G \times_H X) /\!/ G$ if $H$ is a closed subgroup of $G$.
\end{itemize}

With this is mind, we discuss some generalizations of GIT to
non-reductive group actions, including various notions of ``stable"
and ``semistable" points.   A direct approach is to patch together
``nice enough" affine opens which admit finitely generated rings of
invariants.  Another method is to transfer the problem to reductive
GIT by considering the associated reductive $G$-action on
(projective completions of) $G \times_H X$; given a linearized
$H$-action on $X$ in $\PP^n$, we can consider
$$G \times_H X \subseteq G \times_H \PP^n  \cong (G/H) \times
\PP^n, $$ for $G$ a reductive group extending the action of $H$ on
$\PP^n$, and reductive GIT on $G \times_H \PP^n$ or an appropriate
projective completion $\overline{G \times_H \PP^n}$.  In each
approach, and conceptually perhaps this is the crucial point, only
finitely many invariants are really being used to define the
quotients. Somewhat remarkably, when properly formulated these
approaches are compatible, and most of the conditions above are
satisfied --- though the ``enveloping quotient" $X^{ss} \to X/\!/H$
which we obtain fails to be a categorical quotient in general, and
its image may not be a subvariety of $X/\!/H$ but only a dense
constructible subset. As a byproduct we obtain a geometric criterion
for deciding whether or not the ring of invariants $k[X]^H$ is
finitely generated. When $X$ is projective and $k[X]^H$ is finitely
generated then $X/\!/H = \Proj(k[X]^H)$ is a projective variety;
more generally we obtain projective completions $\overline{X/\!/H}$
of $X/\!/H$ which are themselves reductive GIT quotients
$\overline{G \times_H X}/\!/G$ and hence, in principle, amenable to
standard methods for understanding their geometry and topology.

Each affine algebraic group $H$ has a unipotent radical $H_u$, which
by definition is non-trivial precisely when $H$ is non-reductive.
Effectively, one may first quotient by the action of $H_u$, and then
by the induced action of the reductive group $H/H_u$, provided that
the unipotent quotient is sufficiently canonical to inherit an
induced linear action of the reductive group.  So what really must
be understood is unipotent actions, and we shall concentrate on
these.

For convenience, although this is hardly necessary, we work over an
algebraically closed base field $k$ of characteristic $0$, which in
examples will be taken to be $\mathbb{C}$. By ``variety" we mean an
integral scheme of finite type over $k$.  Except where otherwise
stated, we assume $X$ is a projective variety.

This paper is intended to be sufficiently self-contained that
readers with minimal background knowledge will find it accessible,
so it is at times, by necessity or design, somewhat informal. After
this introduction $\S$\ref{sec:Motivation} discusses some problems
in mathematics where non-reductive actions appear: moduli spaces,
exotic affine spaces, Hilbert's fourteenth problem, and so forth.
$\S$\ref{sec:GITBasics} is a brief summary of the idea behind
geometric invariant theory --- using invariants to parametrize
orbits --- followed by an explanation of how key facts about
reductive actions fail for non-reductive ones, and a summary of the
main definitions and results of reductive GIT needed later.
  $\S$4  discusses various \lq intrinsic' ways to
  describe open subsets on which non-reductive
actions admit nice quotients, and compares these to existing
definitions and approaches in the literature.  $\S$5 develops the
approach to non-reductive GIT via reductive GIT on an auxiliary
space; there are three main results. Theorem \ref{thm:main}
summarizes the relationships between the different notions of stable
points, semistable points, and their ``quotients", while Theorem
\ref{thm:main2} provides a stronger conclusion when a particular
criterion is satisfied, which can often be arranged in practice (for
example, when $k[X]^H$ is finitely generated). With these in mind,
Definition \ref{defn:main} sets our notion of stable and semistable
points.  Theorem \ref{thm:fgcriterion} gives a geometric criterion
for deciding when a ring of invariants is finitely generated,
relating it to a (non-effective) stability condition in reductive
GIT. Finally $\S$6 discusses a family of $\CC^+$-actions as a
straightforward example and computes the (intersection) cohomology
of their enveloping quotients $X/\!/H$.

   The authors thank Bob MacPherson for his support and inspiration,
as well as Aravind Asok and Charles Doran for helpful conversations
and comments. Indeed, it was in the first author's time as Bob's
student that this circle of ideas first began to take shape.

\section{Motivation}\label{sec:Motivation}

Reductive group actions are of great significance in algebraic
geometry, but non-reductive actions appear in many important
problems. Over $\CC$ a group $G$ is reductive if and only if it is
the complexification of a maximal compact subgroup $K$, and  many
nice properties  can obtained by exploiting this underlying
compactness. The simplest example here is of course the
complexification $\CC^*$ of the circle $S^1$, and more generally
$GL(n;\CC)$ viewed as the complexification of the unitary group
$U(n)$; in contrast $\CC^+$ has no nontrivial compact subgroups.
Given the ubiquity of translation actions, it is not surprising that
non-reductive groups appear in so many problems, though sometimes
they are well hidden. Two very recent examples include
Bridgeland-Douglas stability conditions in derived categories of
sheaves interpreted, in a special motivating case, as stability for
matrix factorizations using non-reductive actions \cite{HorWalk,
Walk}, and the study of hyperbolic varieties via Griffiths-Green jet
bundles, after Demailly et al.~\cite{Rousseau}; some more classical
examples are given below.

\subsection{Moduli spaces}

In the preface to the first edition of \cite{GIT}, Mumford states
that his goal is ``to construct moduli schemes for various types of
algebraic objects'' and that this problem ``appears to be, in
essence, a special and highly non-trivial case'' of the problem of
constructing orbit spaces for algebraic group actions. More
precisely, when a family of objects with parameter space $S$ has the
local universal property for a given moduli problem, and a group
acts on $S$ such that objects parametrized by points in $S$ are
equivalent if and only if the points lie in the same orbit, then the
construction of a coarse moduli space is equivalent to the
construction of a categorical quotient which is an orbit space for
the action (cf. \cite[Proposition 2.13]{New}). There are many cases
of moduli problems involving non-reductive group actions.  Here are
a few examples.

\begin{itemize}
\item { Moduli of singularities} (or modules over the local ring of a
singularity) as studied in \cite{GPsing1, GPsing2}; here the rough
idea is that translation actions arise from unfoldings of
singularities.

\item Moduli of suitable maps; for example, degree $d$ maps between
projective spaces or with additional constraints (see
\cite[pp.~80-83]{Faunt1} ).

\item Moduli of hypersurfaces in toric varieties.
\end{itemize}

This last example generalizes the very classical moduli problem of
hypersurfaces in projective space $\PP^n$. There the space of
hypersurfaces of degree $d$ is parametrized by
$\PP(Sym^d(\CC^{n+1}))$ and the equivalence is given by the natural
linear action of the reductive group $SL(n+1;\CC)$. There is a
canonical (categorical) projective GIT quotient with an open subset
which is a coarse moduli space for stable hypersurfaces. For a
general toric variety $X$, the analogous parameter space is also a
projective space, but the action is a linear action of a
non-reductive group. Such spaces arise, for example, in the study of
moduli spaces of Calabi-Yau varieties and mirror symmetry
\cite{CoxKatz}. Lacking a theory of non-reductive quotients, the
standard trick is to study an associated ``simplified'' moduli space
which arises as a quotient of a torus action; it is a branched cover
of the actual moduli space and understanding the geometry of its
compactification is complicated and not entirely naturally related
to the actual moduli problem.

\begin{example}
Let $X$ be the weighted projective plane $\PP(1,1,2)$, with
homogeneous coordinates $x,y,z$, and consider the moduli problem of
weighted degree 4 hypersurfaces (which are genus $1$ curves) in this
toric variety $X$. The relevant group action is that of the
automorphism group of $X$, which lifts to a semidirect product of
the unipotent group $(\CC^+)^3$ acting via
$$[x:y:z] \mapsto [x:y:z + \lambda x^2 + \mu xy + \nu y^2] \quad \mbox{
for } (\lambda,\mu,\nu) \in (\CC^+)^3$$ and the reductive group
$GL(2;\CC) \times GL(1; \mathbb{C})$ acting on the $(x,y)$
coordinates and the $z$ coordinate.  A basis for weighted degree $4$
polynomials is
$$\{ x^4, x^3y, x^2y^2,xy^3,y^4,x^2z,xyz,y^2z,z^2 \}.$$
With respect to this basis, the $(\CC^+)^3$ action is linearly
represented as:
$$\left( \begin{array}{ccccccccc}

1 & 0 & 0 & 0 & 0 & \lambda & 0 & 0 & \lambda^2 \\

0 & 1 & 0 & 0 & 0 & \mu & \lambda & 0 & 2\lambda\mu \\

0 & 0 & 1 & 0 & 0 & \nu & \mu & \lambda & 2\lambda\nu + \mu^2 \\

0 & 0 & 0 & 1 & 0 & 0 & \nu & \mu & 2\mu\nu \\

0 & 0 & 0 & 0 & 1 & 0 & 0 & \nu & \nu^2 \\

0 & 0 & 0 & 0 & 0 & 1 & 0 & 0 & 2\lambda \\

0 & 0 & 0 & 0 & 0 & 0 & 1 & 0 & 2\mu \\

0 & 0 & 0 & 0 & 0 & 0 & 0 & 1 & 2\nu \\

0 & 0 & 0 & 0 & 0 & 0 & 0 & 0 & 1

\end{array}  \right)$$
\end{example}

\subsection{Affine geometry}

An abundant supply of translation symmetries means that
non-reductive actions occur everywhere in affine geometry. In
particular the automorphism group ${\rm Aut}(\AAA^n)$ of
$n$-dimensional affine space is a surprisingly rich and mysterious
structure, which ties into many of the famous questions in the area.
To name one example, the {\em Jacobian Conjecture} is equivalent to
the statement that an automorphism $F = (F_1, \ldots, F_n)$ of
$\AAA^n$ naturally induces a set of $n$ locally nilpotent
derivations $\frac{d}{dF_1}, \ldots, \frac{d}{dF_n}$, i.e., a
unipotent action on $\AAA^n$ \cite[$\S$2.2]{vdEssen}. Another
instance is the existence of {\em exotic affine spaces} (varieties
diffeomorphic but not algebraically isomorphic to $\mathbb{A}^n$),
where constructions of examples directly or implicitly make heavy
use of unipotent actions. Perhaps the cleanest to write down is the
Russell cubic three-fold:
$$X = \{ (w,x,y,z) \in \mathbb{C}^4 \ | \ x^2 w + x + y^2 + z^3 = 0 \}, $$
which Makar-Limanov \cite{ML} proves is exotic effectively by
showing it has a smaller set of $\mathbb{C}^+$ actions than
$\mathbb{C}^3$ does (more precisely, the intersection of all the
invariant subrings for the $\mathbb{C}^+$ actions is not simply the
constants). Perhaps more strikingly, arbitrary dimensional families
of non-isomorphic exotic affine spaces can be constructed using
quotients by free unipotent actions \cite{Wink1,AD}.
\begin{ex}\label{ex:q-affquot}
Let $\phi_t: \mathbb{C}^5 \rightarrow \mathbb{C}^5$ be the quadratic
free $\mathbb{C}^+$ action given by:
$$ \phi_t(w_1, w_2, w_3, w_4, w_5) = (w_1, w_2, w_3 + t w_1, w_4 + t w_2,
    w_5 + t(1+w_1 w_4 - w_2 w_3)). $$
    This action on $\mathbb{C}^5$ has a geometric quotient which is
    diffeomorphic to $\mathbb{C}^4$ and quasi-affine but {\em not}
    affine.  The quotient can be thought of as the complement to the cotangent
    space at a point in the cotangent bundle to $S^4$, where the
    bundle is embedded as an affine quadric hypersurface.
\end{ex}

   Moreover Winkelmann \cite{Wink2} has shown that the study of
   the coordinate rings of quasi-affine varieties over any field $k$
   (which are, of course, not necessarily finitely generated as $k$-algebras) is
   precisely the study of the invariant subrings of affine varieties under
   affine algebraic group actions (indeed $k^+$ actions).
   \begin{thm} {\rm (Winkelmann)} If $R$ is an integrally
   closed $k$-algebra then the following are equivalent:
   \begin{itemize}
   \item there exists a quasi-affine irreducible, reduced $k$-variety $X$
such
   that $R \cong k[X]$;
   \item there exists an irreducible, reduced $k$-variety $X$ and a
subgroup
   $G$ of ${\rm Aut}(X)$ such that $R \cong k[X]^G$;
   \item there exists an affine irreducible, reduced $k$-variety $X$ and a
   regular action of $G=k^+$ on $X$ such that $R \cong k[X]^G$.
   \end{itemize}
   \end{thm}

\subsection{Classical invariant theory}

The fourteenth of Hilbert's problems posed at the 1900 International
Congress of Mathematicians was the following question: If an
algebraic group acts linearly on a polynomial ring in finitely many
variables, is the ring of invariants always finitely generated? The
answer is yes for reductive groups (and for some non-reductive
groups, in particular for $\CC^+$), but Nagata \cite{Nagata} showed
that the answer is no, in general, though counterexamples have not
been easy to find. Nagata's original counterexample was an action of
$(\CC^+)^{13}$. Much later Mukai found an action of $(\CC^+)^3$
where the ring of invariants is not finitely generated \cite{Muk},
and made further generalizations in \cite{Mukai, Mukai2}; see
\cite{CT} for some very recent related results on
$(\CC^+)^2$-actions. Popov \cite{Popov} used Nagata's counterexample
to show that if an algebraic group $G$ is not reductive then there
is an affine $G$-variety $X$ such that $k[X]^G$ is not finitely
generated.

\section{GIT Basics}\label{sec:GITBasics}

If $X$ is a normal quasi-projective $G$-variety, then there exists a
$G$-linearization \cite{Dolg}[Theorem 7.3] for the $G$-action on
$X$; that is, a $G$-equivariant embedding in some projective space
together with a lift of the $G$-action to the (homogeneous)
coordinate ring\footnote{Note that if $X$ is not normal then this is
not necessarily true; for example, the nodal cubic curve is a
$\mathbb{C}^*$-equivariant completion of $\mathbb{C}^*$, but not of
any linear $\mathbb{C}^*$ orbit.}. The main goal of geometric
invariant theory is to provide a natural algebraic variety
(depending only on the choice of linearization) which parametrizes
$G$-orbits in an affine or projective variety $X$ by using invariant
functions (or sections) in the (homogeneous) coordinate ring
$A=k[X]$. When the ring $A^G$ of {\em all} invariants is a finitely
generated $k$-algebra, as is the case if $G$ is reductive, we obtain
an affine variety $\Spec(A^G)$ and a morphism
$$F : X=\Spec(A) \longrightarrow \Spec(A^G)$$ for
affine $X$, or
$$F : X= \Proj(A) \dashrightarrow \Proj(A^G)$$ for
projective $X$, induced by the inclusion $A^G \rightarrow A$.  In
the affine case $F$ is a dominant morphism, whereas in the
projective case it is just a dominant rational map; $F$ is not
defined precisely where all the invariants vanish, known as the {\em
unstable set}.  If $X$ is affine and $\{f_i \}_{1 \leq i \leq n}$ is
a finite generating set for $A^G$ , then one may realize $F$ as the
morphism to affine space $(f_1, \ldots, f_n): X \rightarrow
\mathbb{A}^n$; more precisely, this morphism is the composition of F
with an embedding of $\Spec(A^G)$ in $\mathbb{A}^n$. A similar
statement holds in the projective case, except the grading on the
image must be compatible with the degrees of the $f_i$ so that the
image lies in a weighted projective space.

When $A^G$ is {\em not} finitely generated we can still consider
$\Spec(A^G)$ and $\Proj(A^G)$ but only as schemes, not varieties.
\begin{defn}
Let $X$ be a quasi-projective variety, and $L$ an ample line bundle
endowed with a lifting of the $G$-action.  A {\em finite separating
set of invariants} is a collection of invariant sections $\{f_1,
\ldots, f_n \}$ of positive tensor powers of $L$ such that, if $x,y$
are any two points of $X$, and $f$ is an arbitrary invariant section
of $L^{\otimes k}$ for some $k
> 0$, then
$$
( f(x) = f(y), \forall \, f   ) \Leftrightarrow ( f_i(x) =
f_i(y), \forall i = 1, \ldots, n ).
$$
\end{defn}
\begin{rem}
Other authors consider a similar notion.  See \cite[Section
2.3.2]{DerKem} for a comparison.
\end{rem}
Since $X$ is Noetherian one sees that finite separating sets of
invariants always exist (the analogous result in \cite{DerKem} is
Theorem 2.3.15). When $X$ is affine any finite collection of
invariant functions $f_1,\ldots,f_n$ defines a $G$-invariant map
$(f_1,\ldots,f_n):X \to \AAA^n$ (likewise in the projective case),
and one could hope that when $\{f_1,\ldots,f_n\}$ is a finite
separating set of invariants the image of $(f_1,\ldots,f_n)$ might
be a variety independent of the choice of $\{f_1,\ldots,f_n\}$.
Unfortunately this is not true in general, but a refinement of the
idea does work (cf. Proposition \ref{patching2} below).

For a general affine $G$ over a characteristic zero field, there is
a semi-direct product decomposition $G \cong G_r \ltimes G_u$ into a
reductive subgroup $G_r$ and the unipotent radical $G_u$ of $G$, and
the corresponding rings of invariants satisfy $A^G =
(A^{G_u})^{G_r}$. (Note that the characteristic zero hypothesis and
ensuing semi-direct product structure are just conveniences; in any
characteristic $G_u$ is a normal unipotent subgroup of $G$ such that
$G/G_u$ is reductive, and then $A^G = (A^{G_u})^{G/G_u}$). So the
key to non-reductive geometric invariant theory lies in unipotent
actions.

\subsection{Comparison of reductive and non-reductive}

To see how non-reductive actions differ from reductive ones, it
suffices to work with unipotent groups.  Unipotent actions have
several appealing features.
\begin{itemize}
\item Every orbit for a unipotent group action on a
quasi-affine variety is closed \cite{Rosen}, so unipotent actions
are always closed in the sense of Mumford; this is helpful for GIT
since invariants cannot distinguish an orbit from another orbit in
its closure.
\item Every homogeneous space of a unipotent group, and hence every orbit of
a unipotent group action, is isomorphic to an affine space
\cite[Proposition 8.4.1]{KMT}.
\item A connected unipotent group
over a characteristic $0$ field has no proper finite subgroups,
hence a unipotent action has no points with non-trivial but finite
isotropy groups. (This follows from embedding the unipotent group as
a subgroup of upper triangular matrices in some linear
representation.)
\end{itemize}

However there are also several ways in which unipotent actions
behave less well than reductive actions. We list some key properties
of reductive invariant theory which make GIT so effective, and
follow them with unipotent counter-examples.

\begin{prop}\cite[Corollary 1.2]{GIT}\label{prop:sepclosedorbits}
Let $X$ be an affine variety with a reductive $G$-action.  Given two
disjoint $G$-invariant closed subvarieties $Y_1$ and $Y_2$ in $X$,
there exists an invariant $f \in k[X]^G$ such that $f(Y_1) = 0$ and
$f(Y_2) = 1$.
\end{prop}

\begin{prop}\cite[Theorem 1.1(3)]{GIT} Let $X$ be a $G$-invariant affine
subvariety of $\mathbb{A}^n$ where $G$ is reductive and acts on
$\mathbb{A}^n$. Then given $f \in k[X]^G$ there exists an extension
$F$ of $f$ to an invariant on $\mathbb{A}^n$.  More precisely, if
$I_X$ is the defining ideal of $X$, then $k[X]^G =
k[\mathbb{A}^n]^G/I_X \cap k[\mathbb{A}^n]^G$.
\end{prop}

\begin{ex}\label{ex:sym2}
Let $\mathbb{C}^+$ act on $X=Sym^2(V)$ via its inclusion as upper
triangular matrices in the defining two-dimensional representation
$V$ of $SL(2;\mathbb{C})$.\\

\begin{figure}[h!]\label{fig1}
\begin{center}
\includegraphics{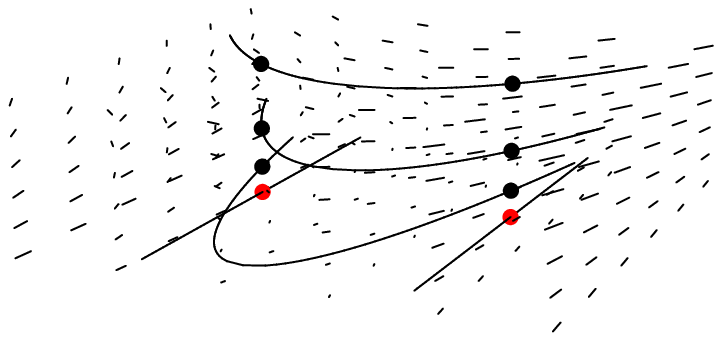}
\caption{$\CC^+$-orbits in $X=Sym^2(V)$}
\end{center}
\end{figure}

Here the invariants are easy to describe with respect to the usual
coordinates $x_0$, $x_1$, $x_2$ on $X=Sym^2(V)$: we have
$A^{\mathbb{C}^+} = \mathbb{C}[x_0, x_1^2 - x_0 x_2]$.  Observe
that, even though this is an action on an affine space, the
invariants do {\em not} separate all closed orbits -- the pairs of
lines $\{ x_0 = 0, x_1 = \pm a \}$ for any given $a \neq 0$ are not
distinguished, nor are any of the fixed points of the line $\{x_0 =
0, x_1 = 0\}$. Furthermore, observe that {\em not all} invariants of
the hyperplane defined by $x_0 = 0$ lift to invariants of $X \cong
\mathbb{C}^3$; for example $x_1$ is an invariant function of the
$\mathbb{C}^+$ action on the hyperplane, but does not extend to an
invariant of the action on $X$.

This example illustrates as well two important principles: firstly,
that non-closed orbits are not the only source of non-separated
orbit spaces; and secondly, that set-theoretic freeness does not
imply scheme-theoretic freeness, as the action of $\mathbb{C}^+$ on
the complement of the fixed point line $\{ x_0 = 0, x_1 = 0 \}$ is
not even separated --- because the parabolic orbits in Figure 1
degenerate into pairs of disjoint orbits in the $(x_0,x_1)$-plane
--- let alone scheme-theoretically free \cite[Definition 0.8]{GIT}
(that is, set-theoretically free and proper \cite[$\S$6.3 Lemma
8]{EdGr}).

\end{ex}

\begin{prop}\cite[Theorem 3.3]{Dolg} \label{prop:Nagata} (Nagata)
Let $G$ be a reductive group acting on an affine variety $X =
\Spec(A)$.  Then $A^G$ is a finitely generated $k$-algebra.
\end{prop}

\begin{ex}[Nagata counter-example]
Consider the $n$-fold direct sum $V^{\oplus n}$ of the defining
representation $V$ for $SL(2;\mathbb{C})$, and look at the subgroup
of $GL(V^{\oplus n})$ corresponding to the direct sum of the $n$
Borel subgroups (each one is a semi-direct product of $\mathbb{C}^+$
and $\mathbb{C}^*$):
$$ \left( \begin{array}{ccccccccc}

c_1 & a_1 & 0 & \ldots & \ldots & \ldots & \ldots & \ldots & 0 \\

0 & c_1 &  0  & \ldots & \ldots & \ldots & \ldots & \ldots & 0 \\

0 &  0  & c_2 & a_2    &   0    & \ldots & \ldots & \ldots & 0 \\

0 &  0  &  0  & c_2    &   0    & \ldots & \ldots & \ldots & 0 \\

\vdots & \vdots & \vdots & \vdots & \vdots & \ddots & \vdots &
\vdots & \vdots \\

0 &  0  & \ldots & \ldots & \ldots & \ldots  & 0 & 0  & c_n
\end{array} \right)
$$
with $a_i \in \mathbb{C}$ and $c_i \in \mathbb{C}^*$. Demand that
$\prod_i c_i = 1$.  Then for appropriate $n$ and a closed subgroup
$G$ cut out by an appropriate set of linear conditions on the $a_i$,
the ring of invariants $\mathbb{C}[x_1, \ldots, x_n, y_1, \ldots,
y_n]^G$ is not finitely generated \cite[Theorem 4.3]{Dolg}.
\end{ex}

\begin{prop} Let $G$ be a reductive group acting on an affine variety $X =
\Spec(A)$, or respectively a projective variety $X = \Proj(A)$. Then
the quotient map $q: \Spec(A) \longrightarrow \Spec(A^G)$, or respectively
   $q: \Proj(A)^{ss} \longrightarrow \Proj(A^G)$,
is surjective.
\end{prop}

\begin{ex} This proposition says that the GIT quotient of an affine
(projective)
variety by a reductive group action is affine (projective), whereas
in contrast Example \ref{ex:q-affquot} provides an affine variety whose
quotient by a non-reductive group action
 is quasi-affine but {\em not} affine.
\end{ex}

\subsection{Geometric and categorical quotients}

Recall (from \cite[Chapter 2, $\S$4]{New}, for example) that when a
group $G$ acts on a variety $X$ a {\em categorical quotient} of $X$
by $G$ is a morphism $\phi:X \to Y$ from $X$ to a variety $Y$ which
is $G$-invariant (that is, constant on $G$-orbits) and has the
property that any other $G$-invariant morphism $\tilde{\phi}: X \to
\tilde{Y}$ factors as $\tilde{\phi} = \chi \circ \phi$ for a unique
morphism $\chi:Y \to \tilde{Y}$. An {\em orbit space} for the action
is a categorical quotient $\phi:X \to Y$ such that $\phi^{-1}(y)$ is
a single $G$-orbit for each $y \in Y$, and a {\em geometric
quotient} is an orbit space $\phi:X \to Y$ with the following {\em
good} properties: it is an affine morphism such that

(i) if $U$ is open in $Y$ then
$$ \phi^*: k[U] \to k[\phi^{-1}(U)] $$
induces an isomorphism of $k[U]$ onto $k[\phi^{-1}(U)]^G$, and

(ii) if $W_1$ and $W_2$ are disjoint closed $G$-invariant
subvarieties of $X$ then their images $\phi(W_1)$ and $\phi(W_2)$ in
$Y$ are disjoint closed subvarieties of $Y$.

Thus an orbit space is a geometric quotient if and only if it is a
good categorical quotient.

\subsubsection{GIT for reductive group actions on projective
varieties}\label{sssec:red-on-proj}

Let $X$ be a projective variety over $k$ and let $G$ be a reductive
group acting on $X$. For Mumford's geometric invariant theory we
also require a {\em linearization} of the action; that is, a line
bundle $L$ on $X$ and a lift\footnote{When there is no risk of
confusion we will use $L$ to denote the linearization as well as the
underlying line bundle.} of the action of $G$ to $L$. When $L$ is
ample we can assume without essential loss of generality that for
some projective embedding $X \subseteq \PP^n$ the action of $G$ on
$X$ extends to an action on $\PP^n$ given by a representation
$$\rho:G\rightarrow GL(n+1;k),$$
where $L$ is the hyperplane line bundle on $\PP^n$. There is an
induced action of $G$ on the homogeneous coordinate ring
$$A = k[X] \overset{def}{=} \bigoplus_{k \geq 0} H^0(X, L^{\otimes k})$$
of $X$, which, when  $X \subseteq \PP^n$ as above, is the quotient
$k[x_0,...,x_n]/ \mathcal{I}_X $ of the polynomial ring
$k[x_0,...,x_n]$ by the ideal $\mathcal{I}_X$ generated by the
homogeneous polynomials vanishing on $X$. The subring $A^G$ of $A$
consisting of the elements of $A$ left invariant by $G$ is a graded
$k$-algebra, which by Nagata's theorem is finitely generated because
$G$ is reductive \cite{Nagata}, so we can define $X/\!/G$ (or
$X/\!/_L G$ when the dependence on the linearization $L$ is to be
made explicit) to be the variety $\Proj (A^G)$ associated to this
ring of invariants $A^G$. The inclusion of $A^G$ in $A$ defines a
rational map $\phi$ from $X$ to $X/\!/G$, but because there may be
points of $X \subseteq \PP_n$ where every nonconstant $G$-invariant
homogeneous polynomial vanishes, this map will not in general be a
morphism.

We define the set $X^{ss}(L)$ (abbreviated to $X^{ss}$ when there is
no risk of confusion) of {\em semistable} points for the action of
$G$ on $X$ with respect to the ample linearization $L$ to consist of
those $x \in X$ for which there exists some $k>0$ and $s \in H^0(X,
L^{\otimes k})^G$ not vanishing at $x$. Then the rational map $\phi$
restricts to a surjective $G$-invariant morphism from the open
subset $X^{ss}$ of $X$ to the projective variety $X/\!/G$, and
$\phi:X^{ss} \to X/\!/G$ is a categorical quotient for the action of
$G$ on $X^{ss}$. Set-theoretically, $X/\!/G$ is the quotient of
$X^{ss}$ by the equivalence relation for which $x$ and $y$ in
$X^{ss}$ are equivalent if and only if the closures
$\overline{O_G(x)}$ and $\overline{O_G(y)}$ of the $G$-orbits of $x$
and $y$ meet in $X^{ss}$.

In order to obtain a geometric quotient, we define a {\em stable}
point for the linear action of $G$ on $X$ to be a point $x$ of
$X^{ss}$ with a neighborhood in $X^{ss}$ such that every $G$-orbit
meeting this neighborhood is closed in $X^{ss}$, and is of maximal
dimension equal to the dimension of $G$ (a ``properly stable point''
in the sense of \cite[Definition 1.8]{GIT}). If $U$ is any
$G$-invariant open subset of the set $X^s= X^s(L)$ of stable points
of $X$, then $\phi(U)$ is an open subset of $X/\!/G$ and the
restriction $\phi|_U :U \to \phi(U)$ of $\phi$ to $U$ is a universal
geometric quotient (it remains a geometric quotient under base
change) for the action of $G$ on $U$. In particular, there is a
geometric quotient $X^s/G = \phi(X^s)$ for the action of $G$ on
$X^s$, and if $X^s$ is nonempty then $X/\!/G$ can be thought of as a
projective completion of the quasi-projective variety $X^s/G$:

$$\begin{array}{ccccc}
     X^s & \subseteq & X^{ss} & \subseteq & X \\
         & {\rm open} &  & {\rm open} & \\
    \downarrow & & \downarrow & & \\
    &  &  &  &  \\
    X^s/G & \subseteq & X^{ss}/\sim & = & X/\!/G. \\
         & {\rm open} & & &
\end{array}$$
\begin{rem}
Each of $X^s, X^{ss}$, and $X/\!/G$ remains unchanged if the
linearization $L$ is replaced by $L^{\otimes m}$ for any $m > 0$.
\end{rem}
\begin{rem} \label{maximality} When $X$ is a connected nonsingular
projective variety and $G$ is a connected reductive group acting on
$X$ with generic stabilizer having dimension 0, then a $G$-invariant
open subset $U$ of $X$ is proper and there exists a quasi-projective
geometric quotient $\phi:U \to Y$ if and only if $U \subseteq
X^s(L)$ for some linearization $L$ of the action of $G$ (see
\cite[Converse 1.13]{GIT}).
\end{rem}

\begin{rem} \label{rem1}
Let $T$ be a maximal torus of a reductive group $G$ acting linearly
on $X$. Then the subsets $X^{ss}$ and $X^s$ of $X$ are characterised
by the following properties (essentially the Hilbert-Mumford
criteria for stability and semistability) which make them easy to
identify:

(i) A point $x \in X$ is semistable (respectively stable) for the
linear action of $G$ on $X$ if and only if $gx$ is semistable
(respectively stable) for the action of $T$ on $X$ for {\em every}
$g \in G$.

(ii) If the maximal torus $T$ of $G$ acts diagonally on $X\subseteq
\PP_n$ with weights $\a_0,...,\a_n$, then a point $x = [x_0,...,x_n]
\in X$ is semistable (respectively stable) for the action of $T$ if
and only if the convex hull
$${\rm Conv} \{ \a_i:x_i \neq 0\}$$
in $\liets$ contains 0 (respectively contains 0 in its interior)
where $\liets$ is the vector space dual of the Lie algebra of $T$.
\end{rem}

\begin{example} \label{projective space}
Consider the linear action of $\CC^*$ on $X=\PP^n$, with respect to
the hyperplane line bundle on $\PP^n$, where the linearization $L_0$
is given by the representation
$$t \mapsto \diag(t^2,t^2,\ldots,t^2,1)$$
of $\CC^*$ in $GL(n+1;\CC)$. The same action of $\CC^*$ on $\PP^n$
has other linearizations with respect to the hyperplane line bundle;
let $L_+$ denote the linearization given by the representation
$$t \mapsto \diag(t^3,t^3,\ldots,t^3,t)$$
of $\CC^*$ in $GL(n+1;\CC)$ and let $L_-$ denote the linearization
given by the representation
$$t \mapsto \diag(t,t,\ldots,t,t^{-1})$$
of $\CC^*$ in $GL(n+1;\CC)$. Then
$$X^{ss}(L_+) = X^s(L_+) = \emptyset,$$
while
$$X^{ss}(L_0) = \CC^n \mbox{ , } X^s(L_0) = \emptyset,$$
and
$$X^{ss}(L_-) = X^s(L_-) = \CC^n \setminus \{ 0 \}.$$
\end{example}

\subsubsection{Reductive GIT and quasi-projective varieties}
\label{sssec:red-qproj}

If a $G$-action on a variety $X$ has a categorical quotient $\phi:X
\to Y$ then its restriction to a $G$-invariant open subset of $X$ is
not necessarily a categorical quotient for the action of $G$ on $U$,
as the following simple example shows.

\begin{example}\label{ex:notfunctorial}
Let the multiplicative group $\CC^*$ of $\CC$ act on $\CC^n$ as
multiplication by scalars. Since the origin lies in the closure of
every orbit, it follows that every $G$-invariant morphism
$\phi:\CC^n \to Y$ is constant and hence that the constant map from
$\CC^n$ to a point is a categorical quotient for the action. However
the restriction of this constant map to $\CC^n \setminus \{ 0 \}$ is
not a categorical quotient, since the natural map $\CC^n \setminus
\{ 0 \} \to \PP^{n-1}$ is a nonconstant $G$-invariant morphism, and
indeed is a categorical quotient for the action of $\CC^*$ on $\CC^n
\setminus \{ 0 \}$. Thus in Example \ref{projective space} above we
have $\PP^n /\!/_{L_-} \CC^* = \PP^{n-1}$ while $\PP^n /\!/_{L_+}
\CC^*$ is empty and $\PP^n /\!/_{L_0} \CC^*$ is a point.
\end{example}

Similarly, although following Mumford \cite{GIT} we can define
stable and semistable points for any linear action of a reductive
group $G$ on a quasi-projective variety $X$, it is not necessarily
the case that $U^{ss} = U \cap X^{ss}$ or that $U^s = U \cap X^s$
whenever $U$ is a $G$-invariant open subset of $X$.

\begin{defn}\label{defn:s/ssred}
Let $X$ be a quasi-projective variety with an action of a reductive
group $G$ and linearization $L$ on $X$. Then $y \in X$ is {\em
semistable} for this linear action if there exists some $m \geq 0$
and $f \in H^0(X, L^{\otimes m})^G$ not vanishing at $y$ such that
the open subset
$$ X_f := \{ x \in X \ | \ f(x) \neq 0 \}$$
is affine, and $y$ is {\em stable} if in addition the action of $G$
on $X_f$ is closed with all stabilizers of dimension 0.
\end{defn}

\begin{rem} \label{rem:GITnotample}
Note that the line bundle $L$ is not required to be ample here. When
$X$ is projective and $L$ is ample and $f \in H^0(X, L^{\otimes
m})^G \setminus \{ 0 \}$ for some $m \geq 0$, then $X_f$ is affine
if and only if $f$ is nonconstant or equivalently $m>0$. Thus
Definition \ref{defn:s/ssred} agrees with
$\S$\ref{sssec:red-on-proj} for projective $X$ and ample $L$.
Moreover when $X \subseteq \AAA^n \subseteq\PP^n$ is {\em affine}
and the linear $G$-action is given by a representation of $G$ in
$GL(n;k)$ embedded in $GL(n+1;k)$ in the usual way by taking a
direct sum with a one-dimensional trivial representation of $G$,
then we have
$$X^{ss} = X$$
and $X/\!/G = \Spec(k[X]^G)$ is a categorical quotient of $X$ (cf.
Theorem \ref{thm 3.2.9} below), while $x \in X$ is stable if and
only if there is some $f \in k[X]^G$ such that the action of $G$ on
$X_f$ is closed with all stabilizers having dimension 0.
\end{rem}

\begin{rem}
One can see from these definitions that three things can go wrong in
relating stable and semistable points for open immersions $U
\subseteq X$: (1) invariants do not necessarily extend, and even if
an invariant $f$ extends then (2) affineness of $X_f$ is neither a
necessary nor a sufficient condition for $U_f$ to be affine, and (3) the
action on $U_f$ may be closed with all stabilizers of dimension 0
without the same being true of $X_f$.
\end{rem}

\begin{thm}[Mumford] \cite[Theorem 1.10]{GIT} \label{thm 3.2.9}
Let $X$ be an algebraic scheme over $k$ with an action of a
reductive group $G$ and linearization $L$ on $X$. Then $X^{ss}$ has
a quasi-projective categorical quotient $\phi:X^{ss} \to X/\!/G$
which restricts to a geometric quotient $\phi:X^s \to \phi(X^s)$ of
$X^s$, where $\phi(X^s)=X^s/G$ is an open subset of $X/\!/G$.
\end{thm}

\begin{rem}
The maximality property of $X^s$ described in Remark
\ref{maximality} extends to nonsingular quasi-projective varieties
$X$.
\end{rem}

\section{Generalizing GIT: Intrinsic viewpoints}

Throughout, unless otherwise stated, $X$ will be a projective
variety endowed with an ample line bundle $L$, on which an affine
algebraic group $H$ acts linearly. As discussed previously, $H$ may
be assumed to be connected unipotent.

\subsection{Global approaches}

This subsection will be fairly informal, serving three purposes --
to set the stage as naively as possible, to present some
definitions, and to give a brief survey of some of the literature.

The inclusion $k[X]^H \hookrightarrow k[X]$ induces a rational
``quotient" map of schemes $q: X \dashrightarrow \Proj(k[X]^H)$. The
image of $q$ is a constructible subset, i.e. a finite union of
locally closed subschemes.  The goal is to understand $q$ well
enough to construct canonically associated ``quotient" varieties.
    There are a couple of interrelated approaches to analyzing
$q$ and possible quotients:
\begin{itemize}
\item  take large enough finitely generated approximations to $k[X]^H$,
and

\item  understand $q$ restricted to open sets on which the action is nice.
\end{itemize}

For simplicity of notation we shall assume that the generic
stabilizer has dimension 0; the discussion in general is not
significantly different.

\subsubsection{Naive stability}

As in reductive GIT for projective varieties with ample line
bundles, we might think of the semistable points as the points where
some $H$-invariant does not vanish; that is, the points in the
domain of definition of $q$.
\begin{defn}
A point $x \in X$ is {\em naively semistable} if $q(x)$ is
well-defined.  The set of naively semistable points $X^{nss}$ is the
domain of definition of $q$.
\end{defn}

Consider an increasing filtration of $A= k[X]^H$ by finitely
generated subrings $A_i$ with $\underset{\rightarrow}{\lim} (A_i) =
A$. The inclusions $A_i \rightarrow A_{i+1} \rightarrow \dots
\rightarrow A \rightarrow k[X]$ induce rational maps
$$q_i: X \dashrightarrow \Proj(A)
\dashrightarrow \dots \dashrightarrow \Proj(A_{i+1}) \dashrightarrow
\Proj(A_i).$$ Because $X$ is Noetherian, the unstable set $X
\setminus X^{nss}$ is cut out by finitely many invariants, so for
large enough $i$, the domain of definition of $q_i$ is $X^{nss}$.
Furthermore, when $i$ is large enough, $A_i$ contains a finite
separating set of invariants; the natural rational map $q_{i+1,i}:
\Proj(k[X]^H_{i+1}) \dashrightarrow \Proj(k[X]^H_i)$ then maps
$q_{i+1}(X^{nss})$ bijectively onto $q_i(X^{nss})$.

It therefore makes sense to look for open subsets $U$ of $X^{nss}$
whose images $q_i(U)$ are varieties independent (up to isomorphism)
of $i$, once $i$ is sufficiently large.  Upper semi-continuity of
$q_i$ for large $i$ determines a finite filtration of $X^{nss}$ by
opens $U_j$ indexed according to the maximal dimension $j \geq
\dim(H)$ of components of the fiber at a point; provided that $A_i$
contains a separating set of invariants this filtration is
independent of $i$.  Restricting one's attention to the smallest of
these open subvarieties, $U_{\dim(H)}$, is a good way to find
analogs of ``stable" points, since $q$ can only be an orbit map (and
hence possibly a geometric quotient) on subsets of $U_{\dim(H)}$.

\begin{defn}
Let $X$ be a projective variety with a linear $H$ action. A point $x
\in X^{nss}$ is said to be {\em almost stable} if its stabilizer
group is trivial and $\dim(q^{-1}(q(x))) = \dim(H)$.  The set of
almost stable points we denote by $X^{as}$.
\end{defn}

\begin{rem}
If $k[X]^H$ is finitely generated and $\Proj(k[X]^H)$ is normal,
then by standard results (see \cite[Proposition 18.4]{Borel}) the
morphism $q$ restricted to $X^{as}$ is an open mapping, so its image
is a variety.  Restricting $q$ further to $U=q^{-1}(V)$ where $V$ is
the union of open subsets of $q(X^{as})$ over which all geometric
fibres are single orbits --- what might be called the ``{\em naively
stable}" set --- it follows \cite[Proposition 6.6]{Borel} that we
obtain a geometric quotient (in fact this holds in general, see
Remark \ref{rem:DR-GP}).
\end{rem}

\begin{rem}
Fauntleroy \cite{Faunt2} considered an analog of this concept of
almost stability: for $Y$ quasi-affine and normal on which a
connected unipotent group $H$ acts with finite stabilizers, and $p:
Y \rightarrow \Spec(k[Y]^H)$ the canonical morphism, if
  $\dim(p^{-1}(p(y)))=\dim(H)$ then he
calls $y$ ``semistable" for the action of $H$.  He shows that $p$
restricted to the open subvariety of $Y$ consisting of these
semistable points is a categorical quotient (in the category of
varieties). Fauntleroy's results applied to open subvarieties of an
affine cone over a normal projective variety $X$ with a linear
$H$-action show that
  $q(X^{as})$ is a
variety.
\end{rem}

Another option is to find ``quotient" varieties $Q$, obtained by
throwing out selected closed $H$-invariant subvarieties of
$\Proj(A_i)$ for sufficiently large $i$, which are as close as
possible to the image of $q_i$ --- ideally while retaining some
other good feature, such as $k[X]^H = k[Q]$. To that end, removing
those (finitely many) codimension 1 divisors in $\Proj(A_i)$ which
contain dense open subsets disjoint from $\im(q_i)$ is a reasonable
approach, although some points of $\im(q_i)$ are typically removed
in the process and some points not in $\im(q_i)$ may remain: this
is, very roughly, what we do later with strong reductive envelopes
(see Definition \ref{defn:strong} and, with more stringent
conditions, Theorem \ref{thm:main2}).  In any such construction,
whether the quotient variety is suitably canonical becomes a serious
issue.

\begin{rem}\label{rem:Winkelss}
In the case when $X$ is affine, Winkelmann \cite{Wink2} showed there
exists a quasi-affine ``quotient" variety $Q$ which admits a {\em
rational} map from $X$ and such that $k[Q]= k[X]^H$.  The idea here,
too, is removal of codimension 1 divisors in $\Spec(A_i)$ for large
enough $i$.
\end{rem}

\subsubsection{Quotients by free actions}

One can also make use of auxiliary properties of a group action
which ensure the existence of a well-behaved quotient; specifically,
recall that in reductive GIT the action on the stable locus is {\em
proper}.  In the case of a connected unipotent group over a
characteristic $0$ field, the proper actions (being also
set-theoretically free) are exactly \cite[$\S$6.3 Lemma 8]{EdGr} the
{\em scheme-theoretically free} actions \cite[Definition 0.8]{GIT}.
   By Artin \cite{Artin} and Koll\'{a}r \cite{Kollar}, proper
actions admit geometric quotients in the category of algebraic
spaces.  So by considering only those open subvarieties on which the
action is proper, coupled with some condition that ensures the
quotient is a variety, one obtains a notion of stable sets, though
finding a canonical choice of stable set is an issue here.
Fauntleroy \cite{Faunt1} defines a notion of {\em properly stable
actions} of connected unipotent groups on quasi-affine normal
varieties over algebraically closed fields, which yield geometric
quotients in the category of varieties. Conversely, sufficiently
well behaved geometric quotients come from properly stable actions.

\begin{defn}
A linear action of a connected unipotent group $H$ on a normal
quasi-affine variety $U$ is {\em properly stable} if the $H$-action
on $U$ is proper and every $x \in U$ has a $H$-invariant open
neighbourhood $W$ such that  the nonempty fibres of the natural map
of schemes $q:W \to \Spec(k[W]^H)$ are $H$-orbits.
\end{defn}

\begin{prop}[\cite{Faunt1}]
Let $U$ be a normal quasi-affine variety over the algebraically
closed field $k$, on which the connected unipotent group $H$ acts.
Assume that the action of $H$ on $U$ is properly stable. Then the
quotient map $q: U \to q(U) \subseteq \Spec(k[U]^H)$ is affine and a
principal $H$-bundle (i.e. a locally trivial geometric quotient by
\cite[Proposition 0.9]{GIT}). Conversely, if $Y$ is a
   variety and $q:U \to Y$ is a principal $H$-bundle,
then the action of $H$ on $U$ is properly stable.
\end{prop}
\begin{rem}
Note that $H$ is a special group, in the sense of Serre and
Grothendieck, so all principal $H$-bundles are in fact Zariski, not
just \'{e}tale, locally trivial \cite[$\S$2.6]{PopVin}.
\end{rem}
\begin{rem}
Unfortunately, a maximal properly stable open subset of an affine
variety $X$ is not in general unique, as Example \ref{bel} below
shows, though when $k[X]$ is a UFD then  being properly stable turns
out to be a local property, so the union of all properly stable open
subsets is the unique maximal properly stable open subset \cite[Thm
2.4]{Faunt1}. By \cite[II Prop 6.2, II Cor 6.16]{Hartshorne} when
$X$ is an affine variety the condition that $k[X]$ be a UFD is the
same as $Pic(X)$ being trivial. When working with a projective
variety $X$ the analogous condition is $Pic(X)=\ZZ$.
\end{rem}

\begin{example} \label{bel}
Let $H=\CC^+$ act on $\PP^3$ via
$$t[x_0:x_1:x_2:x_3] = [x_0: x_1: x_2 + t x_1:x_3 + t(2x_2+ x_0)
+ t^2 x_1]$$ and let $X$ be the smooth irreducible surface
$$x_1 x_3 - x_2^2 - x_0 x_2 = 0.$$
Then $\CC^+$ acts set-theoretically freely but not properly on the
affine variety $W = X \cap \AAA^3$ where $\AAA^3 \subseteq \PP^3$ is
defined by $x_0 \neq 0$, and the (non-separated) quotient is the
affine line $\AAA^1$ with a doubled point; $W$ is the union of two
invariant open subsets on each of which the action is properly
stable with quotient $\AAA^1$ (see \cite{FauntMajid}; cf. \cite{DR}
for a similar example). Here points of an open subset of $X$ on
which the $H$-action is properly stable are not necessarily naively
stable in the sense of Definition \ref{naivelystable}.
\end{example}

\subsection{Gluing local quotients $X_f /\!/ H$}\label{sssec:gluing}

In reductive GIT, the set $X^s$ of $G$-stable points of $X$ is a
union of affine invariant hypersurface complements $X_f$ on which
the action is closed and orbits have maximal dimension, and the
geometric quotient of $X^s$ can be constructed by patching together
the corresponding affine varieties $\Spec (k[X_f]^G)$. Unipotent
actions on affine varieties are always closed, but when a unipotent
group $H$ acts linearly on a projective variety $X$ and $f$ is an
invariant section there is no guarantee that $k[X_f]^H$ is a
finitely generated $k$-algebra, that the natural quotient map $q$
from $X_f$ to $\Spec(k[X_f]^H)$ is surjective, or even that the
image of $q$ is a variety.

However one knows that the {\em field} of invariant rational
functions $k(X)^H$ is finitely generated (by a theorem of
Rosenlicht, see \cite[Theorem 6.2]{Dolg}), so there is an invariant
open, $X_f$, for which $k[X_f]^H$ is finitely generated. When $g$ is
another invariant then $X_f \cap X_g = X_{fg}$ has a finitely
generated ring of invariants $k[X_f]^H[g^{-1}]$. We will see that
taking the union of all such open affines $X_f$ for which the
natural map  $q:X_f \to \Spec(k[X_f]^H)$ has sufficiently good
properties, and patching the associated maps $q$, yields canonical
open sets with nice quotients. Patching works here, in contrast with
Example \ref{bel}, since any orbit in $X_f$ is distinguished from
any orbit in the complement of $X_f$ by the invariant $f$ itself.

Since we would like a stable set to have a geometric quotient, it is
natural to impose the condition that $q:X_f \to \Spec(k[X_f]^H)$
should be a geometric quotient.

\begin{defn} \label{naivelystable} Let $X$ be a projective variety and let
$H$ be a connected unipotent group acting linearly on $X$ with
respect to an ample line bundle $L$. The set of {\em naively stable}
points of $X$ (with respect to the linearization $L$) is
     $$X^{ns} = \bigcup_{f \in I^{ns}} X_f$$ where
$$I^{ns} = \{f
\in H^0(X,L^{\otimes m})^H \mbox{ for some }m>0 \ | \ k[X_f]^H
\mbox{ is finitely generated, and }$$
$$  q: X_f \longrightarrow
\Spec(k[X_f]^H) \mbox{ is a geometric quotient} \}.$$
\end{defn}

\begin{prop} \label{patching}
$q(X^{ns})$ is a quasi-projective variety and $q:X^{ns} \to
q(X^{ns})$ is a geometric quotient.
\end{prop}
\begin{proof}
If $f \in I^{ns}$ then $q(X_f)=\Spec(k[X_f]^H)$ is an affine
variety. By the Noetherian property we can choose $f_0,\ldots,f_r
\in I^{ns}$ such that $X^{ns} = \bigcup_{j=0}^r X_{f_j}$, and
without loss of generality we can assume that there is some $m>0$
such that $f_j \in H^0(X,L^{\otimes m})^H$ for $0 \leq j \leq r$.
For each $j$ by the definition of $I^{ns}$ we can choose finitely
many generators $\{ f_{ij}: 1 \leq i \leq q_j \}$ for
$k[X_{f_j}]^H$, which we can express as $f_{ij} =
g_{ij}/(f_j)^{\ell}$ with $g_{ij} \in H^0(X,L^{\otimes \ell m})^H$
for some large $\ell >0$.

Now let $M= r + \sum_{j=0}^r q_j$ and define $s:X^{ns} \to \PP^M$ to
have homogeneous coordinates given by the sections $f_j^{\ell}$
  and $g_{ij}$ (for $0 \leq j \leq r$ and $1 \leq
i \leq q_j$) of $L^{\otimes \ell m}$. If $y_0,\ldots,y_r$ are the
first $r+1$ homogeneous coordinates on $\PP^M$ which pull back to
$f_0^{\ell}, \ldots, f_r^{\ell}$ on $X^{ns}$, then $s$ maps $X^{ns}$
into $\bigcup_{j=0}^r (\PP^M)_{y_j}$ and for $0 \leq j_0 \leq r$ we
have
$$s^{-1}((\PP^M)_{y_{j_0}}) = X_{f_{j_0}}$$
where
$$s:X_{f_{j_0}} \to (\PP^M)_{y_{j_0}} \cong \AAA^M$$
is the composition of $q:X_{f_{j_0}} \to q(X_{f_{j_0}}) =
\Spec(k[X_{f_{j_0}}]^H)$ with the embedding of
$\Spec(k[X_{f_{j_0}}]^H)$ as a closed subvariety of $\AAA^M$ defined
by the subset
$$\{\frac{f_0^{\ell}}{f_{j_0}^{\ell}}, \ldots, \frac{f_0^{\ell}}
{f_{j_0}^{\ell}} \} \cup \{ \frac{g_{ij}}{f_{j_0}^{\ell}}: 0 \leq j
\leq r, 1 \leq i \leq q_j \} $$ of $k[X_{f_{j_0}}]^H$ which includes
the generators
$$\{ f_{ij} = \frac{g_{ij}}{f_{j_0}^{\ell}}: 0 \leq i \leq q_{j_0} \}$$
of $k[X_{f_{j_0}}]^H$. It follows that $s$ is the composition of
$q:X^{ns} \to q(X^{ns}) = \bigcup_{j=0}^r q(X_{f_j})$ with an
embedding of $q(X^{ns})$ as a locally closed subvariety of $\PP^M$,
and moreover that $q:X^{ns} \to q(X^{ns})$ is a geometric quotient
since $q:X_{f_{j}} \to q(X_{f_j})$ is a geometric quotient for $0
\leq j \leq r$ by the definition of $I^{ns}$.
\end{proof}

We can also consider a more stringent condition, that $q: X_f
\longrightarrow \Spec(k[X_f]^H)$ be a locally trivial geometric
quotient (that is, an $H$-principal bundle with base
$\Spec(k[X_f]^H)$).

\begin{defn}\label{defn:lts}
The set of {\em locally trivial stable} points is $$ X^{lts} =
\bigcup_{f \in I^{lts} } X_f$$ where
\begin{eqnarray*} I^{lts} = \{f
\in H^0(X,L^{\otimes m})^H \mbox{ for some } m>0 \ | \  k[X_f]^H
\mbox{
is finitely generated, and } \\
    q: X_f \longrightarrow \Spec(k[X_f]^H) \mbox{ is a locally trivial
geometric quotient} \}. \end{eqnarray*}
\end{defn}

Then $X^{lts}$ is an open subset of $X^{ns}$, and by patching
together the natural maps  $q:X_f \to \Spec(k[X_f]^H)$ as in
Proposition \ref{patching} we obtain
   a locally trivial geometric quotient $q:X^{lts} \to q(X^{lts})$
where $q(X^{lts}) \cong X^{lts}/H$ is an open subvariety of
$q(X^{ns}) \cong X^{ns}/H$.

\begin{rem}
Note that for a {\em reductive} group acting with trivial isotropy
on an invariant open subvariety $X_f$ local triviality is {\em not}
a more stringent condition; that is, $X^{lts} = X^{ns}$. Indeed the
quotient map is affine, and so by Mumford \cite[Proposition
0.8]{GIT} it is proper, and so here is scheme-theoretically free
\cite[$\S$6.3 Lemma 8]{EdGr}; then by \cite[Proposition 0.9]{GIT} it
is a principal bundle.
\end{rem}

\begin{rem} \label{rem:DR-GP}
Another option is to ignore the finitely generated condition on
$k[X_f]^H$ altogether, and work with $\Spec(k[X]^H)$ and $\Spec
(k[X_f]^H)$ as schemes and the natural map $q:X_f \to
\Spec(k[X_f]^H) \subseteq \Spec(k[X]^H)$ as a map of schemes.  This
is the approach taken by Greuel and Pfister \cite{GP1}, coinciding
with that of Dixmier and Raynaud \cite{DR}. For stability they
require less than Definition \ref{naivelystable}; they demand that
the image $q(X_f)$ should be an open subscheme of $\Spec(k[X]^H)$
with $q: X_f \to q(X_f)$ an open map whose geometric fibres are just
the $H$-orbits in $X_f$. (Strictly speaking their definition is for
affine varieties, or quasi-affine with a choice of affine
embedding). However they show, using \cite[Proposition 2.2.2]{DR},
that the open subscheme $q(X_f)$ of $\Spec(k[X]^H)$ is in fact then
a variety, and furthermore that the restriction of the natural map
of schemes $q: X \longrightarrow \Spec(k[X]^H)$ to their stable set
onto its image in $\Spec(k[X]^H)$ is a geometric quotient in the
category of varieties. It follows that their definition of stability
is essentially equivalent to $X^{ns}$.
\end{rem}

In a similar manner, we can mimic the GIT construction of $X^{ss}$
and its categorical quotient by patching together the affine
varieties $\Spec(k[X_f]^H)$ for all those open subvarieties $X_f$
with finitely generated rings of $H$-invariants. We continue to assume
that $X$ is a projective variety and that $H$ is a connected unipotent
group acting linearly on $X$ with respect to an ample line bundle $L$.
\begin{defn}
We define the {\em finitely generated semistable set} by $$X^{ss,
fg} =  \bigcup_{f \in I^{ss, fg}} X_f$$ where
$$I^{ss,fg} = \{f
\in H^0(X,L^{\otimes m})^H \mbox{ for some }m>0 \ | \ k[X_f]^H
\mbox{ is finitely generated }   \}.$$
\end{defn}
It is not difficult to check using the proof of Proposition
\ref{patching} that the affine varieties $\Spec(k[X_f]^H)$ for $f
\in I^{ss,fg}$ patch to give a quasi-projective variety which is an
open subscheme of the scheme $\Proj(k[X]^H)$; see Proposition
\ref{patching2} below. We introduce some terminology for such
quotients:
\begin{defn}\label{defn:envelopquot}
Let $q: X^{ss, fg} \rightarrow \Proj(k[X]^H)$ be the natural
morphism of schemes.  Then the {\em enveloped quotient} of
$X^{ss,fg}$ is $q: X^{ss, fg} \rightarrow q(X^{ss,fg})$, where
$q(X^{ss,fg})$ is a dense constructible subset of the {\em
enveloping quotient}
$$X /\!/ H = \bigcup_{f \in I^{ss,fg}}
\Spec(k[X_{f}]^H)$$ of $X^{ss, fg}$.
\end{defn}

\begin{rem}
The morphism $q: X^{ss, fg} \rightarrow X/\!/H$ is not in general
surjective, even when $k[X]^H$ is finitely generated, as we will
see in $\S$6. The enveloping quotient $X/\!/H$ is thus not in general
a categorical quotient of $X^{ss,fg}$.
\end{rem}

\begin{prop} \label{patching2} The enveloping quotient
$X/\!/H$ is a quasi-projective variety with an ample line bundle
$L_H \to X/\!/H$ which pulls back to a positive tensor power of $L$
under the natural map $q:X^{ss,fg} \to X/\!/H$.
\end{prop}
\begin{proof}
This follows from a trivial modification of the proof of Proposition
\ref{patching}, replacing $I^{ns}$ with $I^{ss,fg}$, replacing
$q(X^{ns})$ with $X/\!/H$ and replacing $q(X_f)$ with
$\Spec(k[X_f]^H)$ for any invariant $f$ throughout.
\end{proof}

\begin{prop} \label{projenvquot}
If $k[X]^H$ is finitely generated then $X/\!/H$ is the projective
variety $\Proj(k[X]^H)$.
\end{prop}
\begin{proof}
When $k[X]^H$ is finitely generated then so is $k[X_f]^H$ for any $f
\in I = \bigcup_{m>0} H^0(X,L^{\otimes m})^H$, so $I^{ss,fg} = I$
and hence
$$X/\!/H = \bigcup_{f \in I}
\Spec(k[X_{f}]^H) = \Proj(k[X]^H).$$
\end{proof}

\begin{lem} \label{lem:as-in-fg}
If $X$ is normal, then $X^{as} \subset X^{ss,fg}$.
\end{lem}
\begin{proof}
By definition $X^{as} \subseteq X^{nss} = \bigcup_{f \in I} X_f$
where $I = \bigcup_{m>0} H^0(X,L^{\otimes m})^H$. If $f \in I$ then
$X_f$ is affine, so we know by \cite{Faunt2} (where normality is
assumed) that $q(X^{as}\cap X_f)$ is an open subvariety of
$\Spec(k[X_f]^H)$. Thus it is the complement of a closed subscheme
cut out by functions $f_i \in k[X_f]^H$ which we can assume are of
the form $f_i = g_i/f^{\ell}$ for some large $\ell \geq 0$ where
$g_i \in k[X]^H$, and hence is a union of affine schemes
$\Spec(k[X_f]^H)_{f_i} = \Spec(k[X_{g_i f }]^H)$; as subschemes of
the variety $q(X^{as} \cap X_f)$ these are themselves affine
varieties, and hence each $k[X_{g_i f}]^H$ is finitely generated.
But then $X^{as}$ is the union of the corresponding opens $X_{g_i f}
= q^{-1}(\Spec(k[X_{g_i f}]^H)$
   of $X$, and so $X^{as} \subseteq X^{ss,fg}$ as required.
\end{proof}

\section{Generalizing GIT: Induced reductive actions -- from $H$ to $G$}

\subsection{Stability for $G \times_H X$}

A couple of fundamental results on algebraic homogeneous spaces will
be useful (see \cite[Theorem 4.17]{PopVin} and \cite[Corollary
2.8]{Grosshans}, respectively).
\begin{prop}[Matsushima] \label{prop:Matsushima}
Given a reductive $G$, then $G/H$ is affine if and only if $H$ is
reductive.
\end{prop}
\begin{prop} \label{Matsushima2}
If $H$ is unipotent, then $H$ is an {\em observable} subgroup of
$G$; that is, $G/H$ is a quasi-affine variety.
\end{prop}
\begin{cor}
If $H$ is a positive dimensional unipotent group and $G$ is
reductive, then $G/H$ is quasi-affine but not affine.
\end{cor}

Choose a reductive group $G$ with the connected unipotent group $H$
as a closed subgroup.  A linearization of an $H$-action on a
projective (or quasi-projective) $X$ with respect to an embedding of
$X$ in $\PP^n$ gives us a representation of $H$ in $SL(n+1; k)$; if
we assume, as without loss of generality we may, that this
representation is faithful, then we can always choose $G$ to be
$SL(n+1;k)$. Let $G \times_H X$ denote the quotient of $G \times X$
by the free action of $H$ defined by $h(g,x)=(g h^{-1}, hx)$. There
is an induced $G$-action on $G \times_H X$ where $G$ acts on itself
by left multiplication.
\begin{rem} With respect to the
natural $G$-equivariant projection to the homogeneous space $G/H$
given by $[g,x] \mapsto gH$, the algebraic space $G \times_H X$ is a
fiber bundle with fibers isomorphic to $X$ presented as sets of the
form $gX$; it is not merely an algebraic space, but is in fact an
algebraic variety \cite{PopVin}[Theorem 4.19].
\end{rem}
If the action of $H$ on $X$ extends to $G$ we have an isomorphism of
$G$-varieties $G \times_H X \cong (G/H) \times X$ given by
\begin{equation} \label{9Febiso} [g,x] \mapsto (gH, gx). \end{equation}

When $X$ is affine we have
$$k[G \times_H X]^G \cong k[X]^H.$$
When $X$ is projective (or more generally quasi-projective) and $L
\rightarrow X$ is a very ample $H$-linearization inducing an
embedding of $X$ in $\PP^n$, and $G$ is a subgroup of $SL(n+1; k)$,
then we get  a very ample $G$-linearization on $G \times_H X$
by pulling back $\mathcal{O}_{\PP^n}(1)$ via
the sequence:
$$ G \times_H X \hookrightarrow G \times_H \mathbb{P}^n  \cong
(G/H) \times \mathbb{P}^n  \to \mathbb{P}^n,
$$  taking the trivial bundle on the quasi-affine variety
$G/H$; by choosing a $G$-equivariant embedding of $G/H$ in an affine
space $\mathbb{A}^m$ with a linear $G$-action we get a $G$-equivariant
embedding of $G \times_H X$ in
$$\mathbb{A}^m \times \mathbb{P}^n \subset \mathbb{P}^m \times \mathbb{P}^n
\subset \mathbb{P}^{nm+m+n}.$$
 By abuse of notation let us also call this $G$-linearization
$L$.  Then
$$ k[G \times_H X]^G =
  \bigoplus_{m \geq 0} H^0(G \times_H X, L^{\otimes m})^G \cong
\bigoplus_{m \geq 0} H^0(X, L^{\otimes m})^H = k[X]^H.$$

\begin{rem} \label{rem:BorelTransf}
If the action of $H$ on $X$ extends to $G$ we have $G \times_H X
\cong (G/H) \times X$ via the isomorphism (\ref{9Febiso}) above, in
which case this identification of rings of invariants is sometimes
known as the Borel transfer principle (see \cite[Lemma 4.1]{Dolg}).
\end{rem}

Let $i: X \rightarrow G \times_H X$ be the closed immersion given by
$x \mapsto [e,x]$.
\begin{defn}
Let the set of {\em Mumford stable points} for a linear $H$-action
on a quasi-projective variety $X$ be $X^{ms} = i^{-1}((G \times_H
X)^s)$.  Similarly, the set of {\em Mumford semistable points} for
the linear $H$-action on $X$ is $X^{mss} = i^{-1}((G \times_H
X)^{ss})$. Here $(G \times_H X)^s$ and $(G \times_H X)^{ss}$ are
defined as in Definition \ref{defn:s/ssred} for the induced linear
action of $G$ on $G \times_H X$.
\end{defn}

These sets admit geometric and categorical quotients, respectively,
by standard reductive GIT. But because $H$ is unipotent something
quite interesting happens.
\begin{lem}\label{lem:ms/mss}
Every Mumford semistable point is Mumford stable; that is,  $X^{ms}
= X^{mss}$.
\end{lem}
\begin{proof}
Let $O$ be an $H$-orbit in $X^{mss} \setminus X^{ms}$, so that $G
\times_H O$ is a $G$-orbit in an affine hypersurface complement $(G
\times_H X)_F$ which is either closed of non-maximal dimension or
not closed in $(G \times_H X)_F$. In each case there is a unique
closed orbit $G \times_H O'$ in the closure of  $G\times_H O $ in
$(G \times_H X)_F$, which is of non-maximal dimension. In an affine
variety a closed $G$-orbit is of course affine, so by Proposition
\ref{prop:Matsushima} the corresponding stabilizer is a reductive
subgroup of $G$.  Hence any point of $G \times_H O'$ must have a
positive dimensional reductive stabilizer, which is a contradiction
since the stabilizer is a subgroup of $H$ and so is unipotent. Thus
there are no Mumford semistable points which are not Mumford stable;
that is, $X^{ms} = X^{mss}$.
\end{proof}
\begin{rem} \label{rem:q-affNOTaff}
When coupled with the fact that $(G \times_H X)_F$ is always
quasi-affine (see Corollary \ref{cor:q-aff}), it follows from this that
whatever notion one comes up with for ``strictly semistable" orbits
that are not Mumford stable, such orbits must lie in {\em
quasi-affine but not affine} $(G \times_H X)_F$. This is a useful
geometric observation, to be used again in the main theorem below.
\end{rem}

{\em A priori} the definitions of $X^{ms}$ and $X^{mss}$ depend on
the choice of the reductive group $G$, but in fact they do not.
\begin{prop}\label{prop:msintrinsic}
The open subsets $X^{ms}$ and $X^{mss}$ of $X$ are intrinsically
defined, in that they depend only on the linear $H$-action on $X$
and are independent of the choice of $G$.
\end{prop}
\begin{proof}
By Lemma \ref{lem:ms/mss} it is enough to consider  $X^{mss}$.
We may assume that both $H$ and $G$ are subgroups of the
general linear group $GL(n+1;k)$ associated to the projective
embedding of $X$ given by the linearization. Suppose that
$G_1$ and $G_2$ are both subgroups of $G_0 = GL(n+1;k)$ containing $H$.
Let $f$ be an $H$-invariant, and let $F_1,
F_2$ and $F_0$ be the associated $G_1, G_2$ and $G_0$ invariants on $G_1
\times_H X, G_2 \times_H X$, and $G_0 \times_H X$, respectively.  By
definition of $X^{mss}$, we want to show that $(G_1 \times_H
X)_{F_1}$ is affine if and only if $(G_2 \times_H X)_{F_2}$ is
affine.  But then it suffices to show that $(G_1 \times_H X)_{F_1}$ is
affine if and only if $(G_0 \times_H X)_{F_0}$ is affine.

If $(G_1 \times_H X)_{F_1}$ is affine, then because $((G_1 \times_H
X)_{F_1}) \times_{G_1} G_0 = ((G_1 \times_H X) \times_{G_1} G_0)_{F_0} = (G_0
\times_H X)_{F_0}$, and both $G_1$ and $G_0$ are affine, the latter is
affine \cite[p. 196]{PopVin}. Conversely, if $(G_0 \times_H X)_{F_0}$ is
affine, then so is any closed subvariety. But the inclusion of $(G_1
\times_H X)_{F_1}$ in $(G_0 \times_H X)_{F_0}$ is a closed immersion, so
it must be an affine variety.
\end{proof}
The intrinsic nature of $X^{ms}$ and $X^{mss}$ can be seen more
geometrically by the following characterization.
\begin{prop}\label{prop:ms=lts}
The Mumford stable points are precisely the locally trivially stable
points; that is, $X^{ms} = X^{lts}$.
\end{prop}
\begin{proof}
Consider $x \in X^{lts}$.  Then $x \in X_f$ for $f$ an $H$-invariant
section of a positive tensor power of $L$, where $\phi: X_f
\rightarrow X_f/H$ is a principal $H$-bundle and $X_f/H =
\Spec(k[X]^H)_f$ is an affine variety.  By \'{e}tale descent, $\phi$
is a principal $H$-bundle if and only if $\pi:(G \times_H X)_F
\rightarrow X_f/H$ is a principal $G$-bundle, for $F$ the
corresponding $G$-invariant.  But by \cite[Proposition 0.7]{GIT}
this means $\pi$ is an affine morphism, so $(G \times_H X)_F$ is
affine.  It follows that $x \in X^{ms}$.

In the reverse direction we see that $X^{ms} \subseteq X^{lts}$ as
follows. $G$ acts scheme-theoretically freely on $(G \times_H X)^s$
because the action is proper and set-theoretically free
\cite[$\S$6.3 Lemma 8]{EdGr}. Furthermore by reductive GIT this
action has a geometric quotient $(G \times_H X)^s \to (G \times_H
X)^s/G \cong X^{ms}/H$.  Hence by Mumford \cite[Proposition
0.9]{GIT} $(G \times_H X)^s \rightarrow X^{ms}/H$ is a $G$-principal
bundle, so by descent $X^{ms} \rightarrow X^{ms}/H$ is an
$H$-principal bundle, and furthermore it is locally described by
$X_f \rightarrow X_f/H$ for some collection $f$ of $H$-invariants.
\end{proof}
\begin{rem}
It follows that any $H$-orbit in $X^{ns} \setminus X^{lts}$ (see
Definitions \ref{naivelystable} and \ref{defn:lts}) has the property
that the corresponding $G$-orbit lies in a quasi-affine but not
affine $(G \times_H X)_F$ (cf. Remark \ref{rem:q-affNOTaff} above).
\end{rem}
\begin{rem}
When it is convenient to do so we can choose $G$ to be semisimple
throughout this section; indeed we can make a canonical choice of
$G$ as $SL(n+1;k)$, but it is often easier to work with a smaller
semisimple or reductive group $G$.
\end{rem}

\subsection{Completions and reductive envelopes}

The techniques of reductive GIT are most effective when applied to
projective varieties, so it makes sense to choose a $G$-equivariant
projective completion of $G \times_H X$ together with an extension
of the $G$-linearization. However, as discussed in Section
\ref{sssec:red-qproj}, GIT does not behave well with respect to
$G$-equivariant open inclusions.  Recalling from Lemma
\ref{lem:ms/mss} that the $G$-semistable points of $G \times_H X$
are all $G$-stable, we may summarize the issues as follows:
\begin{enumerate}
\item Some invariants may not extend to the boundary, so
$G$-stable points of $G \times_H X$ could become unstable in
$\overline{G \times_H X}$.

\item Even if a given invariant $F$ does extend from $G \times_H X$
to $\overline{G \times_H X}$, the open subvariety $(G \times_H X)_F$
may be affine when $(\overline{G \times_H X})_F$ is not, so again
$G$-stable points in $G \times_H X$ could become unstable in
$\overline{G \times_H X}$.

\item Alternatively, if an invariant $F$ extends from $G \times_H X$
to $\overline{G \times_H X}$, then $(\overline{G \times_H X})_F$
might be affine when $({G \times_H X})_F$ is not, so unstable points
of $G \times_H X$ could become semistable in $\overline{G \times_H
X}$.

\item If the $G$-action on $(G \times_H X)_F$ is closed the
$G$-action on $(\overline{G \times_H X})_F$ need not be, so
$G$-stable points in $G \times_H X$ may become strictly semistable
in $\overline{G \times_H X}$.
\end{enumerate}

\begin{rem}
Recall that if $L \rightarrow Y$ is an ample line bundle, $Y$ is
projective, and $F$ is a section of $L$, then $Y_F$ is affine.  So
the second issue above is not a problem for ample extensions of the
$G$-linearization. For ample extensions, it also follows that the
third issue above can be refined since unstable points of $G
\times_H X$ cannot become stable in $\overline{G \times_H X}$ and
hence can only be strictly semistable or unstable in $\overline{G
\times_H X}$: it follows from Proposition \ref{prop:sepclosedorbits}
that if $O$ is a $G$-orbit in $G \times_H X$ whose points are stable
in $\overline{G \times_H X}$ then there is an invariant $F$ which is
nonzero on $O$ and zero on the boundary, so that $(G \times_H X)_F =
(\overline{G \times_H X})_F$ is affine; thus the points of $O$ are
semistable (and indeed stable by Lemma \ref{lem:ms/mss}) in $G
\times_H X$.
\end{rem}

The first task is to find a completion $\overline{G \times_H X}$
together with an extension of the line bundle $L$, such that
``sufficiently many" invariants $F_i$ extend over the boundary with
$\overline{G \times_H X}_{F_i}$ affine.

\begin{defn} \label{defn:fullysep}
Let $X$ be a quasi-projective variety with an $H$-action with
linearization $L$, and let $G$ be any reductive group containing
$H$. A finite separating set $S$ (in the sense of Definition 3.0.1)
of invariant sections of positive tensor powers of $L$ is a {\em
finite fully separating set of invariants} for the linear $H$-action
on $X$ if

(i) for every $x \in X^{ms}$ there exists $f \in S$ with associated
$G$-invariant $F$ over $G \times_H X$ such that $x \in (G \times_H
X)_{F}$ and $(G \times_H X)_{F}$ is affine; and

(ii) for every $x \in X^{ss,fg}$ there exists $f \in S$ such that $x
\in X_f$ and $S$ is a generating set for $k[X_f]^H$.
\end{defn}

\begin{rem}
The proof of Proposition \ref{prop:msintrinsic} shows that this
definition is independent of the choice of $G$.
\end{rem}
We now make a sequence of definitions to address the issues
enumerated above.

\begin{defn}\label{defn:envelope}
Let $X$ be a quasi-projective variety with a linear $H$-action with
respect to an ample line bundle $L$ on $X$. Fix a closed immersion
$H \hookrightarrow G$ for $G$ a reductive group, and fix a
$G$-equivariant projective completion $\overline{G \times_H X}$ of
$G \times_H X$ together with a $G$-linearization $L'$ which
restricts to the $H$-linearization $L$ on $X$. If every
$H$-invariant $f$ in some finite fully separating set of invariants
$S$ for the $H$-action on $X$ (in the sense of Definition
\ref{defn:fullysep}) extends to a $G$-invariant section of a tensor
power of $L'$ over $\overline{G \times_H X}$, we will call
$\overline{G \times_H X}$ together with the $G$-linearization $L'$ a
{\em reductive envelope} of the linear $H$-action on $X$ (with
respect to the finite fully separating set of invariants $S$).
\end{defn}

\begin{defn}\label{defn:ampleenvelope}
With notation as above, if there exists such an $S$ for which every
$f\in S$ extends to a $G$-invariant section $F$ over $\overline{G
\times_H X}$ such that $(\overline{G \times_H X})_F$ is affine, then
we say that $(\overline{G \times_H X}, L')$ is a {\em fine reductive
envelope}. If $L'$ is ample then we say this is an {\em ample
reductive envelope}.
\end{defn}
\begin{rem} \label{rem:ampleisfine}
Note that an ample reductive envelope is automatically fine.
\end{rem}

Issue (4) above can be overcome when the following stringent
condition holds, from which other good properties will also follow
(see Theorem \ref{thm:main2}).
\begin{defn}\label{defn:strong}
With the notation of Definition \ref{defn:envelope} and given a
reductive envelope $( \overline{G \times_H X}, L')$ with respect to
a finite fully separating set of invariants $S$, let $D_1, \dots,
D_r$ denote the codimension 1 components of the boundary of $G
\times_H X$ in $\overline{G \times_H X}$.  If every $f \in S$
extends to a $G$-invariant $F$ over $\overline{G \times_H X}$ which
vanishes on each component $D_j$, then $( \overline{G \times_H X},
L')$ will be called a {\em strong} reductive envelope of the linear
$H$-action on $X$.
\end{defn}

\begin{prop}\label{prop:envexists}
Given  a faithful $H$-action on a quasi-projective variety $X$ with an ample
linearization $L$, there always exists an integer $m>0$ and an ample reductive envelope
for the induced linearization on $L^{\otimes m}$.
\end{prop}
\begin{proof}
The linearization provides an embedding of $X$ into a projective
space such that $H$ acts as a subgroup of the associated general
linear group; without loss of generality $L$ is the hyperplane line
bundle on this projective space. By the Noetherian property a finite
fully separating set of invariants  $S$ exists.  Moreover there is
some $M>0$ such that there is a finite fully separating set $S$ for
which every $f \in S$ is an invariant section of $L^{\otimes M}$.
Take the union of $S$ and a basis for those sections of
$L^{\otimes{M}}$ which extend over the projective space in which $X$
is embedded. This set of sections defines an embedding of $X$ into
some larger projective space $\mathbb{P}^n$ such that $H$ acts as a
subgroup of $G=GL(n+1;k)$ and moreover every invariant in $S$
extends to an $H$-invariant on $\mathbb{P}^n$. As in Remark
\ref{rem:BorelTransf}, since the action of $H$ on $\mathbb{P}^n$
extends to $G$ we have $G \times_H \mathbb{P}^n \cong (G/H) \times
\mathbb{P}^n$, so that if $\overline{G/H}$ is any $G$-equivariant
projective completion of $G/H$ then we can identify $G \times_H
\mathbb{P}^n$ with an open subvariety of $\overline{G/H} \times
\mathbb{P}^n$ and hence obtain a $G$-equivariant projective
embedding of $G \times_H \mathbb{P}^n$. Every $f \in S$ extends to
an $H$-invariant on $\mathbb{P}^n$ and hence to a $G$-invariant on
$G \times_H \mathbb{P}^n$, so by using the construction above we can
find a modified $G$-equivariant projective embedding of $G \times_H
\mathbb{P}^n$ such that every $f \in S$ extends to a $G$-invariant
section of the hyperplane line bundle on the ambient projective
space $\mathbb{P}^N$. Then the closure of $G \times_H X$ in
$\mathbb{P}^N$ together with the restriction of
$\mathcal{O}_{\mathbb{P}^N}(1)$ to $\overline{G \times_H X}$ is an
ample reductive envelope for some positive power of the
linearization $L$.
\end{proof}

\begin{cor} \label{cor:q-aff}
Let $f$ be an $H$-invariant of $(X, L)$, and $F$ be the
corresponding $G$ invariant over $G \times_H X$.  Then $(G \times_H
X)_F$ is a quasi-affine variety.
\end{cor}
\begin{proof}
Any $F$ may be included in the set of invariants from the proof of
Proposition \ref{prop:envexists}.  Because the result is an ample
reductive envelope, $(\overline{G \times_H X})_F$ is affine.
\end{proof}

\begin{rem} \label{rem:const-extend}
We will see in Remark \ref{rem:gentle-ample-fg} that smooth
varieties $X$ do not always have smooth ample reductive envelopes;
indeed, we cannot necessarily find ample reductive envelopes
satisfying a weaker condition than nonsingularity which we will call
gentleness (see Definition 5.3.7 below). The heart of the matter is
contained in the combination of Proposition 5.3.9 below with the
observation that if a normal affine variety $X$ has a fine strong
reductive envelope $( \overline{G \times_H X}, L')$ with respect to
a finite fully separating set of invariants $S$ which includes a
nonzero constant function on $X$, then $k[X]^H$ is finitely
generated.  To see this, first note that without loss of generality,
one may take a normalization of the reductive envelope.  Let $f\in
S$ be a nonzero constant function on $X$ and $F$ denote its
associated $G$-invariant on $\overline{G \times_H X}$.  Then the
affine variety $(\overline{G \times_H X})_F$ contains $(G \times_H
X)_F$ as a codimension 2 complement, so by normality $G$-invariant
sections over the latter canonically extend to sections over the
former. Because $G$ is reductive, we know the space of $G$-invariant
sections over the affine variety $(\overline{G \times_H X})_F$ is
finitely generated. But $(G \times_H X)_F = G \times_H X$ since $F$
restricts to a nonzero constant function over $X$, and $G$-invariant
sections over $G \times_H X$ correspond exactly to $H$-invariant
sections over $X$.
\end{rem}

\begin{defn}\label{defn:s/ssbar} Let $X$ be a
projective variety with a linear $H$-action and a reductive envelope
$\overline{G \times_H X}$ in the sense of Definition
\ref{defn:envelope}.  Label the inclusions $i: X \hookrightarrow G
\times_H X$ and $j: G \times_H X \hookrightarrow \overline{G
\times_H X}$. Let the {\em completely stable points} of $X$ with
respect to the reductive envelope be the set $$X^{\overline{s}} = (j
\circ i)^{-1}(\overline{G \times_H X}^s).$$ Let the set of {\em
completely semistable points} be $$X^{\overline{ss}} = (j \circ
i)^{-1}(\overline{G \times_H X}^{ss}).$$
\end{defn}

\begin{rem}
It will follow from Theorems \ref{thm:main} and \ref{thm:main2} that
if $\overline{G \times_H X}$ is normal and together with
$G$-linearization $L'$ forms a fine strong reductive envelope, then
$X^{\overline{s}}$ and $X^{\overline{ss}}$ are independent of the
choice of $\overline{G \times_H X}$ and $L'$.
\end{rem}

\subsection{Main Theorem and Corollaries}

As in Definition \ref{defn:s/ssbar} let $\overline{G \times_H X}$
be a reductive envelope for a projective variety $X$ with a linear
$H$-action, and let
$$i: X \hookrightarrow G \times_H X \text{ and }
j: G \times_H X \hookrightarrow \overline{G \times_H X}$$ denote the
natural injections. Throughout
$$\pi: (\overline{G \times_H X})^{ss} \to
\overline{G \times_H X}/\!/G$$
will denote the GIT quotient map for
the reductive envelope, restricted to open subsets (by abuse of
notation) as appropriate.

\begin{thm}\label{thm:main}
Let $X$ be a normal projective variety with a linear $H$-action, for $H$ a
connected unipotent group, and let $(\overline{G \times_H X},L')$ be
any fine reductive envelope as in Definition
\ref{defn:ampleenvelope}. Then
$$
X^{\overline{s}}  \subseteq  X^{lts} = X^{ms} = X^{mss}
\subseteq  X^{ns}  \subseteq X^{as} \subseteq  X^{ss, fg}  \subseteq X^{\overline{ss}} = X^{nss}.
$$
The stable sets $X^{\overline{s}}$,
$X^{lts} = X^{ms} = X^{mss}$ and $X^{ns}$ admit quasi-projective
geometric quotients, given by restrictions of the quotient map $q =
\pi \circ j \circ i$. The quotient map $q = \pi \circ j \circ i$
restricted to the open subvariety $X^{ss,fg}$ is an enveloped
quotient with $q: X^{ss,fg} \rightarrow X /\!/ H$ an enveloping
quotient in the sense of Definition \ref{defn:envelopquot}. The
enveloping quotient
  $X /\!/ H$ is an open
subvariety of $\overline{G \times_H X}/\!/G$ and there is an ample
line bundle $L_H$ on $X/\!/H$ which pulls back to a tensor power
$L^{\otimes k}$ of the line bundle $L$ for some $k>0$, such that for
any $m>0$ the sections of $L_H^{\otimes m}$ on $X/\!/H$ are
precisely the $H$-invariant sections of the tensor power $L^{\otimes
km}$ of $L$ on $X^{ss,fg}$.
\end{thm}
\begin{rem}
The main content of this theorem can be summarized in the following
diagram
$$ \begin{array}{ccccccccccccc}
X^{\overline{s}} & \subseteq & X^{lts} = X^{ms} = X^{mss} &
\subseteq & X^{ns} & \subseteq X^{as} \subseteq & X^{ss, fg} & \subseteq X^{\overline{ss}} = X^{nss}\\
\Big\downarrow &  & \Big\downarrow &  & \Big\downarrow & &
\Big\downarrow & \Big\downarrow \\ X^{\overline{s}}/H & \subseteq &
X^{lts}/H & \subseteq & X^{ns}/H & \subseteq & X/\!/H & \subseteq
\overline{G \times_H X}/\!/G
\end{array}
$$
in which all the inclusions are open and all the
vertical morphisms are restrictions of
$\pi:(\overline{G \times_H X})^{ss} \to \overline{G \times_H X}/\!/G$,
while every vertical morphism except the last is also a restriction of
the map of schemes $q:X^{nss} \to \Proj(k[X]^H)$) induced by
the inclusion of $k[X]^H$ in $k[X]$.
Note that everything here is a quasi-projective variety with the
exception of $\Proj(k[X]^H)$, which is not a variety
unless $k[X]^H$ is finitely generated as a $k$-algebra.
\end{rem}

\begin{proof}
By definition $X^{ss, fg} \subseteq X^{nss}$ and $X^{ns} \subseteq
X^{as}$. Lemma \ref{lem:as-in-fg} shows that $X^{as} \subseteq
X^{ss,fg}$ and $X^{lts} \subseteq X^{ns}$ by Definition 4.2.3. Lemma
\ref{lem:ms/mss} gives $X^{ms} = X^{mss}$ and by Proposition
\ref{prop:ms=lts}, $X^{lts} = X^{ms}$.

Let $x \in X^{\overline{s}}$ and embed $X$ in $\overline{G \times_H
X}$ in the natural way.  Then for some $G$-invariant $F$ restricting
to an $H$-invariant $f$ on $X$ there is an affine subvariety
$(\overline{G \times_H X})_F$ containing $x$ on which the $G$ action
is proper, and which admits a geometric quotient $\pi$ to
$(\overline{G \times_H X})_F /\!/ G$. The boundary in $(\overline{G
\times_H X})_F$ is a $G$-invariant closed set, so the restriction of
$\pi$ to the complement of the boundary (that is, to $(G \times_H
X)_F = G \times_H X_f$) is also a geometric quotient; furthermore
the restricted action is proper because properness is a local
property on the base, and hence, since it is set-theoretically free,
it is in fact scheme-theoretically free \cite[$\S$6.3 Lemma
8]{EdGr}. It follows from \cite[Proposition 0.9]{GIT} that $\pi$
restricted to $G \times_H X_f$ is a principal $G$-bundle, and hence
by \'{e}tale descent the morphism $q: X_f \rightarrow (G \times_H
X_f)/G = X_f / H$ is a principal $H$-bundle.  In particular, $x \in
X^{lts}$.

It is clear that $X^{\overline{ss}} \subseteq X^{nss}$ because any
$G$-invariant section $F$ over $\overline{G \times_H X}$ restricts
to a $G$-invariant section over $G \times_H X$ and hence to an
$H$-invariant section $f$ over $X$. Now we argue the reverse
containment.  By Definition \ref{defn:ampleenvelope}, there exists a
finite fully separating set of $H$-invariants $\{ f_i: 1 \leq i \leq
r \}$ defining $G$-invariant sections $F_i$ which extend over the
reductive envelope $\overline{G \times_H X}$ such that $(\overline{G
\times_H X})_{F_i}$ is affine. Thus $\cup_i X_{f_i} \subseteq
X^{\overline{ss}}$. By Definition \ref{defn:fullysep}, given $x \in
X$, all the invariant sections $f \in k[X]^H$ vanish at $x$ if and
only if $f_i(x) = 0$ for all $i$. Thus $X^{nss} = \bigcup_{i=1}^r
X_{f_i} \subseteq X^{\overline{ss}}$.

$X^{ns}$ has a geometric quotient by Proposition \ref{patching}, and
hence the same is true of any $G$-invariant open subset of $X^{ns}$,
while $q:X^{ss,fg} \to X/\!/H$ is an enveloping quotient by
Proposition \ref{patching2}. It follows from the proof of
Proposition \ref{patching2} together with condition (ii) in the
definition \ref{defn:fullysep} of a finite fully separating set of
invariants that for some $m,M>0$ there exist $H$-invariants $f_0,
\ldots, f_M \in H^0(X,L^{\otimes m})^H$ such that
$$x \mapsto [f_0(x): \ldots : f_M(x)]$$
induces an embedding of $X/\!/H$ as a quasi-projective subvariety of
$\PP^M$, and moreover $f_0, \ldots, f_M$ extend to $G$-invariant
sections $F_0, \ldots, F_M$ of $(L')^{\otimes m}$ over $\overline{G
\times_H X}$. Since $k[\overline{G \times_H X}]^G$ is finitely
generated, by adding restrictions of $G$-invariants over
$\overline{G \times_H X}$ to $f_0,\ldots,f_M$ we can assume that
$F_0, \ldots, F_M$ define an embedding
$$y \mapsto [F_0(y):\ldots : F_M(y)]$$
of $(\overline{G \times_H X})/\!/G = \Proj(k[\overline{G \times_H
X}]^G)$ as a projective subvariety of the same  projective space
$\PP^M$. The latter embedding is the composition of the first
embedding with the natural map $X/\!/H \to (\overline{G \times_H
X})/\!/G $, so the natural map embeds $X/\!/H$ as a subvariety of
$(\overline{G \times_H X})/\!/G $.

Finally, the existence of an ample line bundle $L_H$ on $X/\!/H$
with the desired properties is proved in Proposition
\ref{patching2}. The identification between invariant sections of
$L$ over $X^{ss, fg}$ and sections of $L_H$ over $X /\!/ H$ follows
trivially by restriction and patching.
\end{proof}

\begin{rem} Recall that a fine reductive envelope always exists
for any linear $H$-action on a projective variety $X$ by Proposition
\ref{prop:envexists} and Remark \ref{rem:ampleisfine}. A choice of
fine reductive envelope $\overline{G \times_H X}$ provides a
projective completion
$$\overline{X/\!/H} = \overline{G \times_H X}/\!/G$$
of the enveloping quotient $X/\!/H$, which in general depends on the
choice of reductive envelope. Note however that when $k[X]^H$ is
finitely generated then $X/\!/H = \Proj(k[X]^H)$ is itself
projective by Proposition \ref{projenvquot} and hence
$${X/\!/H} = \overline{G \times_H X}/\!/G$$
for any fine reductive envelope $\overline{G \times_H X}$.
\end{rem}

\begin{rem}
The Hilbert-Mumford criterion from reductive GIT (see Remark 3.2.5)
can be used (at least when $L'$ is ample, but see also
\cite{Hausen})
to determine $X^{\overline{s}}$ and $X^{\overline{ss}}$
and to analyze the orbit structure of $X^{\overline{ss}} \setminus
X^{\overline{s}}$ and the quotient map $q: X^{\overline{ss}} \to
X/\!/H$ (cf. \cite{KirwanWeinstein}), which leads in examples to
identification of the intermediate stable and semistable sets
$X^{ms} = X^{mss} = X^{lts} \subseteq X^{ns} \subseteq X^{as}
\subseteq X^{ss, fg}$. In especially good situations we have
$X^{\overline{ss}} =X^{\overline{s}}$ (as in the example in the
final section $\S$6 below with $n$ odd) and then
$$X^{\overline{s}} = X^{ms} = X^{mss} =  X^{lts} = X^{ns} =
X^{as} = X^{ss, fg} = X^{\overline{ss}}.$$
\end{rem}

\begin{thm}\label{thm:main2}
If $\overline{G \times_H X}$ is normal and together with a line
bundle $L'$ provides a fine strong reductive envelope for the linear
$H$-action on $X$, then $X^{\overline{s}} = X^{ms}$ and $X^{ss,fg} =
X^{nss}$.
\end{thm}
\begin{proof}
    We show that $X^{ms} \subseteq X^{\overline{s}}$. Consider $x \in
X^{ms}
= X^{mss}$ and embed $X$ in $G \times_H X$ in the natural way. Then
there exists an $F \in H^0(G \times_H X, L^{\otimes m})^G$ for some
$m>0$ such that $F(x) \neq 0$ and $(G \times_H X)_{F}$ is affine and
the action of $G$ on $(G \times_H X)_{F}$ is closed with all
stabilizers of dimension 0. The section $F$ extends to a
$G$-invariant section over the reductive envelope of the $H$-action
on $X$ and $F$ vanishes on every codimension $1$ component $D_j$ of
the boundary.  We claim that $(\overline{G \times_H X})_{F} = (G
\times_H X)_{F}$: they are both affine, so they are determined by
their coordinate rings, and they differ at most in codimension 2, so
by normality their coordinate rings agree; therefore since one is
contained in the other, they must in fact be equal. Thus $x \in
X^{\overline{s}}$.

We prove now that $X^{nss} \subseteq X^{ss,fg}$.  Consider a finite
fully separating set $S = \{ f_i : 1 \leq i \leq r\}$ associated to
the reductive envelope, and label the corresponding $G$-invariant
sections of powers of $L'$ over $\overline{G \times_H X}$ by $F_i$.
Then $x
  \in (\overline{G \times_H X})_{F_i}$ for one of these
$G$-invariants $F_i$. Since $G$ is reductive and $(\overline{G
\times_H X})_{F_i}$ is an affine variety, it has a finitely
generated ring of invariants. Furthermore, by the definition of a
strong reductive envelope, each $F_i$ vanishes on the codimension
$1$ components $D_j$ of the boundary of $G \times_H X$ in
$\overline{G \times_H X}$. It follows that the complement of the
quasi-affine $(G \times_H X)_{F_i}$ in the affine variety
$\overline{G \times_H X}_{F_i}$ has codimension $2$. Thus by
normality we have $k[X_f]^H = k[(G \times_H X)_F] = k[(\overline{G
\times_H X})_F]^G$ which is finitely generated,  so $x \in X^{ss,
fg}$.  The result now follows from Theorem \ref{thm:main}.
\end{proof}

We will see in $\S$\ref{sssec:strong} that fairly mild conditions on
the singularities of a completion $\overline{G \times_H X}$ are
sufficient to prove the existence of a fine strong reductive
envelope.

\begin{rem}
Let $X$ be an affine variety; then we can define $q: X^{ss,fg} \to
X/\!/H \subseteq \Spec(k[X]^H)$ by analogy with the definition of
$q: X^{ss, fg} \rightarrow X/\!/H \subseteq \Proj(k[X]^H)$ when $X$
is projective, and $X/\!/H$ is a quasi-affine variety with
coordinate ring $k[X/\!/H] \cong k[X^{ss,fg}]^H$. When we have an
ample strong reductive envelope which is normal then by Theorem
\ref{thm:main2}
   $X^{ss,fg} = X^{nss}$ where $X^{nss} = X$ since
$X$ is affine (cf. Remark 3.2.12), so
  $X/\!/H$ is an explicit construction of
Winkelmann's quasi-affine ``quotient" variety \cite{Wink2}, and the
quotient morphism $q:X \longrightarrow X/\!/H$ plays the r\^{o}le of
the morphism (rather than just a rational map) Winkelmann hoped for.
\end{rem}

\begin{defn}\label{defn:main} Let $X$ be a projective variety
equipped with a linear $H$-action. We say that a point $x \in X$ is
{\em stable} if $x \in X^{ms}$ and that it is {\em semistable} if $x
\in X^{ss, fg}$.
\end{defn}

\subsubsection{Constructing strong reductive
envelopes}\label{sssec:strong}

Although Proposition \ref{prop:envexists} guarantees the existence
of an ample (and hence fine) reductive envelope, its method of
construction requires us to identify a finite fully separating set
of invariants, and it is not very clear how to do this in practice.
  The difficulty in building a fine reductive envelope lies in
the trade-offs between the choice of completion, finding a line
bundle so that enough sections extend, and guaranteeing that the
$F=0$ complements are affine. In this subsection we discuss how, for
a fixed sufficiently nice (for example, smooth) completion, one can
find a line bundle so that the sections extend; when this line
bundle is ample, we get a strong ample reductive envelope.  The
technique, with slight adjustment, will be applied in
$\S$\ref{sssec:HtoG} and $\S$\ref{sssec:fg}.

\begin{defn}
Let $Y$ be a normal quasi-projective variety.  We say that a
completion $\overline{Y}$ is {\em gentle} if it is normal and some
integral multiple of each boundary Weil divisor is Cartier.
\end{defn}

\begin{rem}
Over a characteristic $0$ field, as we have assumed from the start,
a gentle completion can always be arranged by equivariant resolution
of singularities.
\end{rem}

Given a gentle $G$-equivariant completion $\overline{G \times_H X}$
of $G \times_H X$, some positive tensor power of any line bundle on
$G \times_H X$ will extend non-canonically across the boundary. One
can check the appropriate cocycle condition \cite[proof of Converse
1.13, page 41]{GIT} to verify that in fact the $G$-linearization
extends as well. If the boundary is codimension at least $2$, then
by normality all invariant sections extend uniquely to invariant
sections over the whole of $\overline{G \times_H X}$.  If the
boundary has codimension $1$ components, the basic idea is to weight
them heavily enough that any given invariant section extends
uniquely over the boundary, vanishing on the codimension 1
components.

To make this precise, let $\{D_j : 1 \leq j \leq r \}$ denote the
collection of codimension $1$ components of the boundary of $G
\times_H X$ in $\overline{G \times_H X}$, let $L$ be the
$G$-linearization on $G \times_H X$ discussed in $S$5.1, and let
$L^{\prime}$ be some chosen extension to a $G$-linearization over a
projective completion $\overline{G \times_H X}$.  Denote by
$L^{\prime}_N$ the induced $G$-linearization on $L^{\prime}[N
\sum_{j=1}^r D_j]$ when $N$ is such that $ND_j$ is Cartier for $1
\leq j \leq r$. Then we have the following proposition.

\begin{prop}\label{prop:crucial}
Let $\overline{G \times_H X}$ be a gentle $G$-equivariant completion
of $G \times_H X$, with a $G$-linearization $L'$ of the $G$-action
which extends the given linearization $L$. Given a finite fully
separating set $S$ of invariants  on $X$, then $(\overline{G
\times_H X}, L^{\prime}_{N})$ is a strong reductive envelope with
respect to $S$ for suitable sufficiently large $N$. If moreover
$(\overline{G \times_H X}, L^{\prime})$ is a fine reductive envelope
with respect to $S$ then $(\overline{G \times_H X}, L^{\prime}_{N})$
is a fine strong reductive envelope with respect to $S$, and hence
Theorems \ref{thm:main} and \ref{thm:main2} apply.
\end{prop}
\begin{proof} (Cf. \cite[proof of Converse 1.13, page 41]{GIT}).
For any given $f \in S$, there exists an $N_f$ such that $f$ extends
to a section of $L^{\prime}_{N_f}$ over the codimension 1 components
$D_j$ of the boundary of $G\times_H X$ in $\overline{G \times_H X}$,
and thus by normality extends to an invariant section $F$ over the
whole of $\overline{G \times_H X}$. We can identify $H^0(\overline{G
\times_H X}, L^{\prime}_n) $ with a subspace of $ H^0( \overline{G
\times_H X}, L^{\prime}_{n+1})$ for any $n$, so that if $f$ extends
to a section $F$  of $L^{\prime}_n$ then it vanishes on each $D_j$
as a section of $L^{\prime}_{n+1}$. Thus taking $N > {N}_f$ forces
$F$ to vanish on each of the codimension 1 components $D_j$ of the
boundary. Consequently, since $S$ is finite, simply take $N > max_{f
\in S} ({N}_f)$; then $(\overline{G \times_H X}, L^{\prime}_{N})$ is
a strong reductive envelope for the linear $H$-action on $X$.

For the final part, observe that the complement of a Cartier divisor
in an affine variety is affine, because any line bundle on an affine
variety is ample \cite[Example II.7.4.2]{Hartshorne}. Thus the
complement of the $D_j$ in each affine $(\overline{G \times_H
X})_{F}$ is affine, and the result follows.
\end{proof}

\begin{rem} \label{rem:gentle-ample-fg}
Unfortunately the following three properties do not in general hold
simultaneously, although we can arrange any two of them: the
completion $\overline{G \times_H X}$ being gentle, the line bundle
$L'$ being ample, and the completion together with the line bundle
forming a reductive envelope (that is, enough invariants extending).
More concretely, if $X$ is affine and there exists a gentle ample
reductive envelope, then $k[X]^H$ must be finitely generated, since
by Remark \ref{rem:const-extend}, it suffices to show there exists a
fine strong reductive envelope with respect to a finite fully
separating set $S$ of invariants which includes a nonzero constant
function over $X$.  Because the completion is gentle, we can
construct $L'_N$ as in Proposition \ref{prop:crucial} above for any
such $S$, and for large enough $N$ it defines a fine strong
reductive envelope.

In particular, the Nagata counter-examples admit no gentle ample
reductive envelopes.

Observe that essentially the same argument shows that if $G/H$ has a
gentle $G$-equivariant affine completion
$\overline{G/H}^{\mathit{aff}}$ then $G/H$ can be represented as a
codimension 2 complement in an affine $G$-variety, and hence $k[G/H]
= k[G]^H$ is finitely generated (cf. \cite[$\S$4]{Grosshans}).
\end{rem}

\begin{rem} \label{amplechoice} Suppose that $L'$ is chosen to be ample
in Proposition \ref{prop:crucial}. If the line bundles $L'_N$ are
also ample for large $N$ (for example, if the divisors given by
positive integral multiples of the $D_j$ are ample), then we can
compare the $G$-linearizations $(\overline{G \times_H X},L'_N)$ for
different $N \geq 0$ using the theory of variation of GIT quotients
\cite{Thaddeus, DolgHu}. Since GIT is unaffected when we replace a
line bundle by any positive tensor power of itself we can consider
$L^{\prime}_N$ for all positive rational values of $N$. We know
\cite{Thaddeus, DolgHu} that the interval $(0,\infty)$ can be
divided into (rational) subintervals $I_j$ such that when $N$ lies
in the interior of a subinterval $I_j$ the GIT quotient (and the
stable and semistable sets) defined with respect to the
linearization $L^{\prime}_N$ depends only on the subinterval in
which $N$ lies. Moreover if $N$ lies on the boundary between two
subintervals then there is a nonempty set of points which are
semistable with respect to the linearization $L'_N$ (for this
particular $N$) but not semistable for the whole family, having as
stabilizer a reductive subgroup $R$ of dimension at least one in
$G$. Any maximal torus of $R$ acts with zero weights on the fibres
of $L'_N$ at all semistable fixed points of $R$. Up to conjugacy
only finitely many subgroups can occur as stabilizers, and their
fixed point sets have only finitely many connected components. From
this and the requirement of zero weights on the fibres of
$L^{\prime}_N$, it follows that only finitely many positive $N$ can
occur as the boundary between two subintervals $I_j$, and hence that
there are only finitely many such subintervals. \end{rem} This
proves:
\begin{lem}\label{lem:s-ssbar/Nindep} For $N \geq 0$ let
$(\overline{G \times_H X},L'_N)$ be an ample strong reductive
envelope for a linear $H$-action on a projective variety $X$. Then
the GIT quotient $\overline{G \times_H X}/\!/G$ and the stable and
semistable sets in $\overline{G \times_H X}$ defined with respect to
the $G$-linearizations $L^{\prime}_N$ are independent of $N$ when
$N$ is sufficiently large.
\end{lem}

\subsubsection{Important case: $H$-action extends to
$G$-action}\label{sssec:HtoG}

If a linear $H$-action on a projective variety $X$ extends to a
linear reductive $G$-action on $X$, then since $G \times_H X \cong
(G/H) \times X$ (cf. Remark \ref{rem:BorelTransf}), the problem of
constructing an ample reductive envelope reduces to understanding
completions of $G/H$ and extensions of the trivial line bundle over
these completions.  Since $G/H$ is quasi-affine, the natural choice
is to take a $G$-equivariant affine closure
$\overline{G/H}^{\mathit{aff}}$, and then projectively complete with
the hypersurface at infinity.  Many practical applications of
non-reductive GIT, for example to moduli spaces of hypersurfaces in
toric varieties, fall into this setting. Two scenarios work out
especially well:
\begin{itemize}

\item if a $G$-equivariant affine closure
$\overline{G/H}^{\mathit{aff}}$ contains $G/H$ as a codimension 2
complement, and

\item if a $G$-equivariant affine closure
$\overline{G/H}^{\mathit{aff}}$ is nonsingular.
\end{itemize}
These are both special cases of a third:
\begin{itemize}

\item if a $G$-equivariant affine closure
$\overline{G/H}^{\mathit{aff}}$
is gentle (that is, given a Weil divisor in the boundary, some
integer multiple is Cartier),
\end{itemize}
though in fact this apparent extra generality is spurious, as we
observed in Remark \ref{rem:gentle-ample-fg}: the existence of a
gentle $G$-equivariant affine closure implies that $G/H$ can be
embedded as a codimension 2 complement in an affine $G$-variety,
which in turn implies that $k[G/H]=k[G]^H$ is finitely generated
\cite[$\S$4]{Grosshans}, and by the Borel transfer principle that
$$k[X]^H \cong k[(G/H) \times X]^G$$
is finitely generated whenever the linear action of $H$ on $X$
extends to a linear $G$-action.

Note that the affine closure $\overline{G/H}^{\mathit{aff}}$ of
$G/H$ can always be chosen to be normal, so we assume that
throughout.  We will see that any of the three conditions above
implies the existence of a strong ample reductive envelope when the
$H$-action extends to a linear $G$-action, and hence, by Theorem
\ref{thm:main2}, that $X^{\overline{s}} = X^{ms}$ and $X^{ss, fg} =
X^{\overline{ss}} = X^{nss}$ for such an envelope, while $X/\!/H =
\Proj(k[X]^H)$ by Proposition \ref{projenvquot}.

Let $\overline{G/H}$ denote the normalized projective completion
$\Proj(k[\overline{G/H}^{\mathit{aff}}] \otimes k[x])$. Let $L'$
denote the $G$-linearization obtained by taking the hyperplane line
bundle on $\overline{G/H}$ (with its natural $G$-linearization)
tensored with the given $G$-linearized ample line bundle $L$ on $X$.
Let $D_{\infty}$ represent the hypersurface at infinity in
$\overline{G/H}$. Note that $D_\infty$ is a Cartier divisor on
$\overline{G/H}$ (corresponding to the hyperplane line bundle), so
$\overline{G/H}$ is a gentle completion of $G/H$ if and only if
$\overline{G/H}^{\mathit{aff}}$ is gentle.

\begin{lem} \label{lem:HtoG-strongenvexists}
  If the given linear $H$-action on $X$ extends to a linear
$G$-action, and if $G/H$ has a gentle $G$-equivariant affine closure
$\overline{G/H}^{\mathit{aff}}$, then $(\overline{G/H} \times
X,L'_N)$ is a strong ample reductive envelope for suitable large
$N$, where $L'_N$ is defined as in Proposition \ref{prop:crucial}.
\end{lem}
\begin{proof}
For  $L'_N$ as in Proposition \ref{prop:crucial} to be defined
requires that $ND_j$ should be a Cartier divisor for each
codimension 1 component $D_j$ of the boundary of $G \times_H X$ in
$\overline{G \times_H X}$. Since any line bundle on an affine
variety is ample, and since the hypersurface at infinity is an ample
divisor, $L'_N$ is ample for every such $N$. The result now follows
immediately from Proposition \ref{prop:crucial}.
\end{proof}

\begin{prop} \label{prop:HtoG-nice}
If the linear $H$-action on a projective variety $X$ extends to a
linear $G$-action, and if $\, \overline{G/H}^{\mathit{aff}}$ is
gentle, then $X^{ms}=X^{\bar{s}}$ and $X^{nss}=X^{\overline{ss}}=
X^{ss,fg}$ has a canonical enveloping quotient $X/\!/H =
\Proj(k[X]^H)$ which is realized as a reductive GIT quotient.
\end{prop}
\begin{proof} This follows from Lemma \ref{lem:HtoG-strongenvexists},
Proposition \ref{projenvquot} and Theorems \ref{thm:main} and
\ref{thm:main2}.
\end{proof}

\begin{rem}
An even more special case is when $H$ is a maximal unipotent subgroup
of $G$; this has been studied by many authors including \cite{Kraft1,Kraft}
and recently \cite{GuilleminSjamaar,Hu}.
\end{rem}

Even when the linear $H$-action fails to extend to a linear
$G$-action on $X$, this is true for the ambient projective space
$\PP^n$, so we obtain the following corollary.

\begin{cor} If
$\overline{G/H}^{\mathit{aff}}$ is gentle, then
$$X^{\overline{s}} = X \cap (\mathbb{P}^n)^{ms}$$
is independent of the choice of $\overline{G/H}^{\mathit{aff}}$.
\end{cor}

\begin{proof}
Since $(G/H) \times X$ is closed in $(G/H) \times \PP^n$, by
reductive GIT \cite[Theorem 1.19]{GIT} we have $X^{\bar{s}} = X \cap
(\PP^n)^{\bar{s}}$, and by Proposition \ref{prop:HtoG-nice} applied
to $\PP^n$ we have $(\PP^n)^{\bar{s}}=(\PP^n)^{ms}$.
\end{proof}

\subsubsection{Consequences for finite generation}\label{sssec:fg}
Let us now assume that $X$ is nonsingular, so that by resolution of
singularities we can find a nonsingular (and hence gentle)
projective completion $\overline{G \times_H X}$. The construction of
$L'_N$ in Proposition \ref{prop:crucial} gives us a necessary and
sufficient condition for the ring of invariants $k[X]^H$ to be
finitely generated; when $L'_N$ is ample, this condition is
 that for large enough $N$ the
codimension $1$ components of the boundary of $G \times_H X$ in
$\overline{G \times_H X}$ are all {\em unstable}.

\begin{thm}\label{thm:fgcriterion}
Let $X$ be a nonsingular projective variety acted on by a connected
unipotent group $H$, and let $L$ be a linearization of the
$H$-action with respect to a projective embedding $X \subseteq
\PP^n$ of $X$. Let $H \hookrightarrow G$ be a closed immersion into
a reductive subgroup of $SL(n+1;k)$. Let $L$ be the induced
$G$-linearization over $G \times_H X$ and let $L'$ be an extension
of $L$ over a gentle completion $(\overline{G \times_H X})$.

Then the ring of invariants $k[X]^H$ is finitely generated if and
only if there exists $ N$ such that, for all $ N^{\prime} > N$ for
which $L'_{N'}$ is defined, any $G$-invariant section of a positive
tensor power of $L^{\prime}_{N^{\prime}}$ vanishes on every
codimension 1 component $D_j$ in the boundary of $G \times_H X$ in
$\overline{G \times_H X}$, where $L'_N = L'[N \sum_j D_j]$.
\end{thm}
\begin{proof}
By restriction $k[\overline{G \times_H X}]^G \subseteq k[G \times_H
X]^G = k[X]^H$.  The forward direction is then a consequence of the
construction of Proposition \ref{prop:crucial} as follows. For large
enough $N$, any given invariant section over $G \times_H X$ extends
and vanishes on each $D_j$. So for large enough $N$ the finitely
many generators of the ring of invariants will all vanish on every
$D_j$, hence every element of $k[\overline{G \times_H X}]^G$
vanishes on every $D_j$.

The reverse direction follows by proving that for any such $N$ the
ring of invariants $k[X]^H \cong k[G \times_H X]^G$ is isomorphic to
the ring of invariants $k[\overline{G \times_H X}]^G$, defined with
respect to the linearization $L'_N$, which is finitely generated
since $\overline{G \times_H X}$ is a projective variety acted on
linearly by the reductive group $G$. This isomorphism arises since
any invariant section $s$ over $G \times_H X$ of $L'_{N}$ extends as
in the proof of Proposition \ref{prop:crucial} above to an invariant
section of $L'_{N'}$ over $\overline{G \times_H X}$ for some $N' >
N$. By hypothesis this section vanishes on each $D_j$ and hence
defines a section of $L'_{N'-1}$ extending $s$. Repeating this
argument enough times we find that $s$ extends to a section of
$L'_N$. The same argument applies to any invariant section $s$ over
$G \times_H X$ of a positive tensor power $(L'_{N})^{\otimes m}$ of
$L'_{N}$, so we have $k[G \times_H X]^G \cong k[\overline{G \times_H
X}]^G$ as required.
\end{proof}

\begin{cor}
In the setting of Theorem \ref{thm:fgcriterion}, if the $L'_N$ are
ample for all $N' > N$, then $k[X]^H$ is finitely generated if and
only if every $D_i$ is unstable for all such $N'$.
\end{cor}
\begin{proof}
For an ample bundle, the complement of the zero set of a section is
affine.  In particular, given an ample $G$-linearization, the set on
which all invariant sections vanish is precisely the unstable set.
\end{proof}

\begin{rem}
For $H$-actions extending to $G$-actions on a nonsingular projective
variety $X$, this necessary and sufficient condition can be made
effective by an explicit construction of a
suitable projective completion for
$G/H$ as in $\S$ \ref{sssec:HtoG}, together with a careful analysis
   of Hilbert-Mumford in this setting (cf. Remark 3.2.5),
 at least when the bundles $L'_N$ are ample. When the
bundles $L'_N$ are not ample the analysis of stability is less
straightforward (cf. \cite{Hausen}). Unfortunately ampleness does
not always occur, as can be seen when $H=\mathbb{C}^+$ and
$G=SL(2;\mathbb{C})$ and $\overline{G/H}$ is the blow-up of $\PP^2$
at $0$.
\end{rem}

\section{Example: $n$ unordered points on $\mathbb{P}^1$}

Let $H=\CC^+$, identified with the group of upper triangular
matrices of the form
$$ \left( \begin{array}{cc} 1 & b\\ 0 & 1 \end{array} \right)$$
in $GL(2;\CC)$, act linearly on $X=\PP^n = \PP(Sym^n(\CC^{2}))$ via
the standard representation of $GL(2;\CC)$ on $Sym^n(\CC^{2})$.

Let $G=SL(2;\CC)$; then we can identify $G/H$ with $\CC^2 \setminus
\{ 0 \}$ via the usual transitive action of $SL(2;\CC)$ on $\CC^2
\setminus \{ 0 \}$ which extends to a linear action on its
projective completion $\PP^2= \overline{G/H}$ with the point
$[1:1:0]$ representing the identity coset $H$. Since the linear
action of $H$ on $X$ extends to $G$ we have $G \times_H X \cong G/H
\times X$, and we are in the setting of $\S$\ref{sssec:HtoG}. Since
$\mathbb{P}^2$ is smooth, this is a gentle completion. Let $L$
denote the hyperplane line bundle on $X=\PP^n$ and let $L_2$ denote
the hyperplane line bundle on $\overline{G/H} = \PP^2$. For any
positive integers $p$ and $q$ there is then an induced linearization
of the action of $G$ on $\overline{G/H} \times X$ with respect to
the line bundle $L^{\otimes p} \otimes L_2^{\otimes q}$. Note this
is a line bundle of the type denoted by $L'_q$ in
$\S$\ref{sssec:HtoG}. In particular, by Lemma
\ref{lem:HtoG-strongenvexists} this provides a strong ample
reductive envelope for large enough $q$.  By Theorem \ref{thm:main2}
this means $X^{\overline{s}} = X^{ms}$ and $X^{nss} = X^{ss, fg}
\overset{\text{def}}{=} X^{ss} = X^{\overline{ss}}$.

We know that $\CC[X]^H$ is finitely generated because $G/H$ is a
codimension 2 complement in its affine closure $\mathbb{A}^2$
\cite[$\S$4]{Grosshans}. We can also see this by applying the
finitely generated criterion, Theorem \ref{thm:fgcriterion}: here
there is only one boundary divisor, namely the product of the line
at infinity with $X$; by the Hilbert-Mumford numerical criterion
(see Remark 3.2.3) it is easy to see that this is unstable, and in
particular that all the invariants vanish there, for sufficiently
large $q$.

By the Hilbert-Mumford numerical criterion, a point of
$\PP(Sym^n(\CC^{2})) \times \{[1:1:0]\} \subseteq
\PP(Sym^n(\CC^{2})) \times \PP^2$, represented by an unordered
sequence $p_1,\ldots,p_n$ of points in $\PP^1$, is stable for this
linear action of $G$ provided that
\begin{itemize}
\item strictly fewer than $\frac{n}{2} + \frac{q}{p}$ of the points
$p_1,\ldots,p_n$
coincide anywhere on $\PP^1$, and
\item strictly fewer than $\frac{n}{2}$ of the points $p_1,\ldots,p_n$
coincide at $[1:0]$;
\end{itemize}
it is $G$-semistable unless
\begin{itemize}
\item strictly more than $\frac{n}{2} + \frac{q}{p}$ of the points
$p_1,\ldots,p_n$ coincide anywhere on $\PP^1$, or
\item strictly more than $\frac{n}{2}$ of the points $p_1,\ldots,p_n$
coincide at $[1:0]$.
\end{itemize}
When $q$ is large compared with $p$ then the first condition is
vacuous in each case, and so a point of $\PP(Sym^n(\CC^{2}))$
represented by an unordered sequence $p_1,\ldots,p_n$ of points in
$\PP^1$ is in $X^{\overline{s}} = X^{s}$  provided that strictly fewer than
${n}/{2}$ of the points coincide at $[1:0]$, and is in
$X^{\overline{ss}} = X^{ss}$ unless strictly more than ${n}/{2}$
coincide there.

Thus when $n$ is odd semistability and stability coincide and we
have a geometric quotient $X^s/H = X^{ss}/H$ which is an open subset
of $X /\!/H = \Proj(k[X]^H) \cong (\overline{G/H} \times X) /\!/G$;
its complement can be identified with the reductive quotient
$X/\!/G$. Observe that a point of $X = \PP(Sym^n(\CC^{2}))$
represented by an unordered sequence $p_1,\ldots,p_n$ of points in
$\PP^1$ is $G$-stable provided that strictly fewer than ${n}/{2}$ of the points
coincide anywhere on $\PP^1$, and is $G$-semistable unless strictly
more than ${n}/{2}$ coincide anywhere on $\PP^1$.

When $n$ is even we have a geometric quotient $X^s/H$ which is an
open subvariety of $X/\!/H \cong (\overline{G/H} \times X) /\!/G$,
again with complement $X/\!/G$, but now the image of $X^{ss}$ in
$X/\!/H$ is not a subvariety: it is the union of the open subvariety
$X^s/H$ and the closed subvariety -- which is in fact a point --
$(X/\!/G) \setminus (X^{s,G}/G)$ where $X^{s,G}$ is the stable set
for the action of $G$ on $X$.

When $n$ is odd $X^s = X^{ns} = X^{as} = X^{ss}$.  Indeed, whenever
$n > 2$ then $X^s = X^{ns} = X^{as}$: by continuity, given the
$\CC^+$-action on $\PP^1$ no $\CC^+$-invariant can distinguish among
orbits consisting of points which correspond to configurations with
${n}/{2}$ points at $[1:0]$, and if $n>2$ this set is more than
1-dimensional (the dimension of $\CC^+$).

Note that the algorithm based on flattening stratifications
described in \cite{GP1} produces a set of stable points much smaller
than $X^{s}$; for these $\CC^+$-actions the algorithm removes the
hyperplane $x_0 = 0$, where $x_0$ is the unique $\CC^+$-invariant
coordinate function.

\subsection{$n = 3$ and $4$}

Let $x_0, \ldots , x_n$ be the usual coordinates on
$Sym^n(\CC^{2})$, so that if $[x_0,\ldots,x_n] \in
\PP(Sym^n(\CC^{2}))$ corresponds to an unordered sequence
$p_1,\ldots,p_n$ of points in $\PP^1$ then $x_0, \ldots, x_n$ are
the coefficients of a homogeneous polynomial of degree $n$ in two
variables whose roots are $p_1,\ldots,p_n$.

When $n$ is small, using Gr\"{o}bner basis techniques (in
particular, an adaptation of an algorithm in \cite[Chapter
1]{vdEssen}) we can explicitly find generators for the ring of
$H$-invariants in the polynomial ring
$$ \bigoplus_{j \geq 0} H^0(X,L^{\otimes j}) = \CC[x_0: \ldots : x_n].$$
For $n=3$ we have four generators
$$x_0,\quad x_0 x_2- x_1^2, \quad x_1^3 + \frac{1}{2}x_0^2 x_3 - \frac{3}{2}x_0 x_1 x_2$$
and $$ 2x_0 x_1 x_2 x_3 - \frac{1}{3}x_0^2 x_3^2 + x_1^2 x_2^2
-\frac{4}{3} x_1^3 x_3 -\frac{4}{3} x_0 x_2^3$$ for the invariants,
so all the invariants vanish if and only if
$$x_0 = x_1 = 0,$$
or equivalently if and only if at least two of the three points
$p_1,p_2,p_3$ coincide at $[1:0]$, as expected from the analysis
above. The enveloping quotient $X/\!/H = \Proj(k[X]^H)$ is a degree
6 hypersurface in the weighted projective space $\PP(1,2,3,4)$ with
equation
$$3 y_0^2 y_3 = 4 y_1^3 - 4 y_2^2$$
and $X^s/H = X^{ss}/H$ is the complement of the point $[0:0:0:1]$ in
this hypersurface.

When $n=4$ we have five generators
$$x_0, \quad x_0 x_2- x_1^2, \quad
x_2^2 + \frac{1}{3} x_0 x_4 - \frac{4}{3} x_1 x_3,
$$
$$x_1^3 + \frac{1}{2}x_0^2 x_3 - \frac{3}{2}x_0 x_1 x_2, \quad $$
and $$x_2^3 - 2 x_1x_2x_3 + x_0x_3^2 + x_1^2x_4 - x_0x_2x_4$$
   so all the invariants vanish if and only if
$$x_0 = x_1 = x_2 = 0,$$
or equivalently if and only if at least three of the two points
$p_1,p_2,p_3,p_4$ coincide at $[1:0]$, as expected. In this case the
boundary $(X/\!/H) \setminus (X^s/H)$  of $X^s/H$ in $X/\!/H$ is
$\PP^1$ and the image of $X^{ss}$ in $X/\!/H$ is the union of the
open subset $X^s/H$ and the point $\infty$ in its complement
$\PP^1$.  The enveloping quotient $X//H = \Proj(k[X]^H)$ is a degree
6 hypersurface in the weighted projective space
$\mathbb{P}(1,2,2,3,3)$ with equation
$$4y_3^2 - 4y_1^3 - y_0^3y_4 + 3y_0^2y_1y_2 = 0.$$

\begin{rem}
For $n > 4$ the problem of computing invariants rapidly becomes
unmanageable.  Indeed already for $n=5$ a computation in Singular
produces $15$ invariant generators and is unable to verify that it
is a complete set.  The difficulty is not surprising, as this
problem is closely related to the well-known one of producing a
generating set of $SL_2$-invariants for binary forms of degree
$n+2$, which is very hard for any $n+2 > 6$.
\end{rem}

\subsection{Cohomology of quotients}

We can study the rational intersection cohomology \cite{IHT,IHT2} of
any reductive GIT quotient $Y/\!/G$, where $Y$ is a projective
variety, by stratifying $Y$ and using equivariant (intersection)
cohomology \cite{JK,JKKW,Kiem,Kir1,K2,K3,K4}. These methods apply in
particular
  to any ample reductive envelope $Y=\overline{G \times_H X}$ of
a linear $H$-action on a projective variety $X$ over $k = \CC$,
allowing us to investigate the intersection cohomology of the
projective completion $Y/\!/G = \overline{X/\!/H}$ of the enveloping
quotient $X/\!/H$; when, as in the examples we are considering,
$k[X]^H$ is a finitely generated $k$-algebra so that $X/\!/H =
\Proj(k[X]^H)$ is projective, this is the intersection cohomology of
the enveloping quotient $X/\!/H$ itself.

The simplest situation is when $Y = \overline{G \times_H X}$ is
nonsingular; this holds in our examples since $Y = \PP^2 \times
\PP^n$. Then we have a $G$-equivariantly perfect stratification $\{
S_{\beta}: \beta \in \mathcal{B} \}$ of $Y$ by $G$-invariant locally
closed nonsingular subvarieties of $Y$ with $S_0 = Y^{ss}$ as an
open stratum, so that
\begin{equation} \label{poincare}
P_t^G(Y^{ss}) = P_t^G(Y) - \sum_{\beta \neq 0} t^{2d_{\beta}}
P_t^G(S_{\beta})
\end{equation}
where $P_t^G(Z) = \sum_{i \geq 0} t^i \dim H^i_G(Z;\QQ)$ is the
$G$-equivariant Poincar\'{e} series of a $G$-space $Z$ and $d_\beta$
is the (complex) codimension of $S_{\beta}$ in $Y$ \cite{Kir1}. When
$G=SL(2;\CC)$ and $Y=\PP^2 \times \PP^n$ then the stratification $\{
S_{\beta} : \beta \in \mathcal{B} \}$ is given by
$$S_0 = Y^{ss}$$
and for $n/2 < j \leq n$
$$S_{j,0} = \{ (w,x) \in \PP^2 \times \PP^n : \mbox{ $x$ represents
$n$ points of $\PP^1$} $$
$$ \mbox{exactly $j$ of which coincide at $[a:b]$ where $y=[1:ta:tb]$
for some $t \in \CC$} \} $$ with codimension $j$ in $Y$, and for $0
\leq j \leq n$
$$S_{j,1} = \{ (w,x) \in \PP^2 \times \PP^n : \mbox{ $x$ represents
$n$ points of $\PP^1$} $$
$$ \mbox{exactly $j$ of which coincide at $[a:b]$ where $y=[0:a:b]$} \} $$
with codimension $j+1$ in $Y$. Equation (\ref{poincare}) gives us
$$P^G_t(Y^{ss}) = (1+t^2 + t^4)(1+t^2 + \cdots + t^{2n})(1-t^4)^{-1}$$
$$ - \sum_{\frac{n}{2} < j \leq n} t^{2j}(1-t^2)^{-1}
- \sum_{0 \leq j \leq n} t^{2j+2}(1-t^2)^{-1}.$$ When $n$ is odd we
have $Y^{ss} = Y^s$ and so $Y/\!/G = Y^{ss}/G = Y^s/G$ and the
Poincar\'{e} polynomial $P_t(X/\!/H) = \sum_{i\geq 0} t^i \dim
H^i(X/\!/H;\QQ)$ of the enveloping quotient $X/\!/H = Y/\!/G$ is
given by
$$P_t(X/\!/H)  =  P^G_t(Y^{ss})  =
[(1+t^2+t^4)(1+t^4+t^8+ \cdots + t^{2n-2}) -(t^{n+1} + t^{n+3} +
\cdots + t^{2n})$$ $$ - (t^2+t^4+\cdots+ t^{2n+2})](1-t^2)^{-1} $$
$$ = (1+t^4 + t^8 + \cdots + t^{2n-2} - t^{n+1} - t^{n+3} - \cdots -
t^{2n}) (1-t^2)^{-1} $$ $$ = 1 +t^2 +2t^4 + 2t^6 + \cdots + [1+
\frac{\min(j,n-1-j)}{2}]t^{2j} + \cdots + 2t^{2n-8} + 2t^{2n-6} +
t^{2n-4} + t^{2n-2} $$ where $[\ \ ]$ denotes the integer part. Thus
$[1+(\min(j,n-1-j))/{2}]$ is the $j$th Betti number of $X/\!/H$ when
$n$ is odd; this is also the $j$th intersection Betti number since
$X/\!/H$ is an orbifold and so its rational intersection cohomology
is the same as its ordinary rational cohomology. The Poincar\'{e}
polynomial of the geometric quotient $X^s/H = X^{ss}/H = (X/\!/H)
\setminus (X/\!/G)$ is given by
$$P_t(X^s/H) = P_t(X/\!/H) - t^4 P_t(X/\!/G) = 1+t^2 + t^4 + \cdots
+ t^{n-1}.$$ We can in fact see this directly, since $X^s$ retracts
onto $\PP^{(n-1)/2}$ and $H$ is contractible.

When $n$ is even then $Y^{ss} \neq Y^s$ and $X/\!/H=Y/\!/G$ is not
an orbifold, which means that there is more work to do to compute
the intersection Betti numbers of $X/\!/H$. We first find the
ordinary Betti numbers of a \lq partial desingularization'
$\widetilde{X/\!/H}$ of $X/\!/H$ (see \cite{K2}) obtained by blowing
up the image in $X/\!/H = Y/\!/G$ of the subvariety $Z$ of $Y$
consisting of those $(w,x) \in \PP^2 \times \PP^n$ with $w=[1:0:0]$
(the origin in $\CC^2 \subseteq \PP^2$) and $x$ representing $n$
points on $\PP^1$, exactly half of which coincide at some $u \in
\PP^1$ and the remaining half coincide at some $v \neq u$. This
partial desingularization $\widetilde{X/\!/H}$ (which has only
orbifold singularities) is itself a projective completion of the
geometric quotient $X^s/H$. It can be represented as
$\widetilde{X/\!/H} = \tilde{Y} /\!/G$ where $\tilde{Y}^{ss}=
\tilde{Y}^s$ is obtained from $Y^{ss}$ by blowing up along its
intersection with $Z$ and removing the proper transform of the
subvariety consisting of those $(w,x) \in \PP^2 \times \PP^n$ with
$x$ representing $n$ points on $\PP^1$ exactly half of which
coincide at some $u =[a:b] \in \PP^1$, and with $w=[1:ta:tb]$ for
some $t \in \CC$. We obtain
$$P_t(\widetilde{X/\!/H}) = P^G_t(\tilde{Y}^{ss})$$
$$ = P_t^G(Y^{ss}) + (t^2 + t^4 + \cdots + t^{2n-2})(1-t^4)^{-1}
- t^n (1+t^2 + \cdots + t^n)(1-t^2)^{-1}$$
$$= 1+2t^2 + 3t^4 + \cdots + (n/2)t^{n-2} + (n/2)t^n + \cdots + t^{2n-2}.$$
Finally from this, using the decomposition theorem of \cite{BBDG},
we can obtain the intersection Poincar\'{e} polynomial $IP_t(X/\!/H)
= \sum_{i \geq 0} t^i \dim IH(X/\!/H;\QQ)$ as
$$IP_t(X/\!/H) = P_t(\widetilde{X/\!/H}) - (t^2+t^4 + 2t^6 + \cdots
+ [\frac{n}{4}]t^{n-2} + [\frac{n}{4}]t^n + \cdots + t^{2n-6} +
t^{2n-4})$$
$$= 1+t^2 + 2t^4 + 2t^6 + \cdots + [(n+2)/4]t^{n-2} + [(n+2)/4]t^n
+ \cdots + 2t^{2n-6} + t^{2n-4} + t^{2n-2}.$$ For more details see
\cite{K3}. The Poincar\'{e} polynomial of the geometric quotient
$X^s/H$ when $n$ is even is given by
$$P_t(X^s/H) = P_t(\widetilde{X/\!/H}) - t^4 P_t(\widetilde{X}/\!/G) =
1+t^2 + t^4 + \cdots
+ t^{n-2}$$ as we can also see directly, since $X^s$ retracts onto
$\PP^{(n-2)/2}$ and $H$ is contractible.

We can also compute intersection pairings in $IH(X/\!/H;\QQ)$ and
$H^*(\widetilde{X/\!/H};\QQ)$, and the ring structure of
$H^*(\widetilde{X/\!/H};\QQ)$ and of $H^*(X/\!/H;\QQ)$ when $n$ is
odd, using the methods of \cite{JK,JKKW}.

\begin{rem} One can work similarly with other (generalized)
cohomology theories, like $K$-theory, or motivic cohomology in the
sense of Voevodsky \cite{ADK1}, at least when $n$ is odd.
\end{rem}

\end{document}